\documentclass[12pt]{article}
\usepackage{amsmath}
\usepackage{amssymb}
\usepackage{eufrak}

\title{On the Constructive Inverse Problem in Differential Galois Theory}
\author{William J. Cook (1), Claude Mitschi (2)\ and  Michael F.
Singer (3)\footnote{Correspondence: Michael F. Singer, North Carolina State University, Department of Mathematics, Box 8205, Raleigh, North Carolina 27695-8205, USA, singer@math.ncsu.edu}}

\begin{document}
\def\QED{\hbox{\hskip 1pt \vrule width4pt height 6pt depth 1.5pt \hskip 1pt}}
\newcounter{defcount}[section]
\setlength{
\parskip}{1ex}
\newtheorem{thm}{Theorem}[section]
\newtheorem{lem}[thm]{Lemma}
\newtheorem{cor}[thm]{Corollary}
\newtheorem{prop}[thm]{Proposition}
\newtheorem{defin}[thm]{Definition}
\newtheorem{remark}[thm]{Remark}
\newtheorem{remarks}[thm]{Remarks}
\newtheorem{ex}[thm]{Example}
\def\GL{{\rm GL}}
\def\SL{{\rm SL}}
\def\Sp{{\rm Sp}}
\def\sl{{\EuFrak {sl}}}
\def\sp{{\EuFrak {sp}}}
\def\gl{{\EuFrak {gl}}}
\def\g{{\EuFrak {g}}}
\def\h{{\EuFrak {h}}}
\def\PSL{{\rm PSL}}
\def\SO{{\rm SO}}
\def\Gal{{\rm Gal}}
\def\Sym{{\rm Sym}}
\def\Gbar{\overline{G}}
\def\gbar{\overline{g}}
\def\Abar{\overline{A}}
\def\Abarbar{\tilde{\tilde{A}}}
\def\calGbar{\overline{\calG}}
\def\calAbar{{\overline{\calA}}}
\def\tildea{\tilde a}
\def\tildeb{\tilde b}
\def\tildeg{\tilde g}
\def\tildee{\tilde e}
\def\tildeA{\tilde A}
\def\e{\epsilon}
\def\calG{{\cal G}}
\def\calH{{\cal H}}
\def\calh{{\cal h}}
\def\calT{{\cal T}}
\def\calC{{\cal C}}
\def\calP{{\cal P}}
\def\calS{{\cal S}}
\def\gS{{\EuFrak S}}
\def\curve{{\rm {\bf C}}}
\def\P1{{\rm {\bf P}}^1}
\def\Ad{{\rm Ad}}
\def\ad{{\rm ad}}
\def\Spin{{\rm Spin}}
\def\Aut{{\rm Aut}}
\def\Int{{\rm Int}}
\def\Al{{\rm A}_{\ell}}
\def\Bl{{\rm B}_{\ell}}
\def\Cl{{\rm C}_{\ell}}
\def\Dl{{\rm D}_{\ell}}
\def\esi{{\rm E}_{6}}
\def\ese{{\rm E}_{7}}
\def\obar{\overline{\omega}}
\def\cf{{\it {cf.\ }}}
\def\ie{{\it {i.e.\ }}}
\def\semi{\ltimes}
\newenvironment{prf}[1]{\trivlist
\item[\hskip \labelsep{\bf #1.\hspace*{.3em}}]}{~\hspace{\fill}~$\square$\endtrivlist}
\newenvironment{proof}{
\begin{prf}{Proof}}{
\end{prf}}
 \def\square{\QED}
\newenvironment{sketchproof}{
\begin{prf}{Sketch of Proof}}{
\end{prf}}

\date{}
\maketitle

\vskip 30pt
\centerline{(1) and (3): Department of Mathematics, North Carolina State University}
\centerline{Raleigh,  USA}

\medskip
\centerline{(2) : Institut de Recherche Math\'ematique Avanc\'ee,
Universit\'e Louis Pasteur}
\centerline{Strasbourg, France}

\bigskip
\begin{abstract}{We give  sufficient conditions for a differential equation to have a given
semisimple group as its Galois group.  For any group $G $ with $G^0
= 
G_1 \cdot\ldots \cdot G_r$ where each $G_i$ is a simple group  of type $\Al$,  $\Cl$,  $\Dl$,  ${\rm E}_6$ or ${\rm E}_7$, 
we construct a differential equation over $C(x)$ having Galois group $G$. }
\end{abstract}

\bigskip
\noindent{\it Keywords:} Linear differential equations, Picard-Vessiot, inverse problem, Galois theory.

\bigskip
\noindent{\it 2000 Mathematics Subject Classification:} 34M50 (primary), 12H05, 12H20 (secondary).

\section{Introduction} 

 In (Singer, 1993)  a large class of linear algebraic groups, including  all groups with semisimple identity component, are shown to occur as Galois groups of differential equations
$\frac{dY}{dx} = AY$ with $A$ an $n\times n$ matrix with coefficients in  
$C(x)$ where $C$ is an algebraically closed field of characteristic zero. 
The proof of this depended heavily on the analytic solution of the Riemann-Hilbert Problem and did
not directly give a way of constructing such an equation\footnote{Once one knows that a group $G$
occurs as a differential Galois group over $C(x)$ for $C$ a recursive algebraically closed field of
characteristic $0$, one can produce an equation with Galois group $G$ by listing all equations and using the
algorithm from Hrushovski (2002) to calculate an equation's Galois group   to test until one is
found. 
In this paper, we are interested in more direct methods.}.  Techniques for constructing an equation
with a given group have been produced for connected solvable groups
by Kovacic (1969, 1971) and for  connected groups in general 
by Mitschi and Singer (1996). For groups that are not connected, Mitschi and Singer (2002) showed 
that one could construct equations having  any solvable-by-finite group as Galois group assuming one
could produce (algebraic) equations having any finite group as Galois group. 
Julia Hartmann (2002) has shown that any linear algebraic group can be realized as the Galois group of a 
linear differential equation over  
$C(x)$ and this proof shows that equations can be constructed once one knows how to construct
equations for reductive groups (the proof uses the results of Singer (1993) in this case as
well).\\
[0.2in]
In this paper we will give a criterion, Proposition~\ref{prop2}, for a differential equation to have
a given semisimple algebraic group as a Galois group.  We will use this proposition to show how one
can construct differential equations with Galois group $G$, where  $G = H \ \semi G^0, \ H$ 
$\mbox{finite and } G^0 = 
G_1\cdot \ldots \cdot G_r$ with each $G_i$ of type $\Al$,  $\Cl$,  $\Dl$,  $\esi$ or $\ese$. 
  The rest of the paper is organized as follows. In Section 2, we discuss the criteria given by Mitschi and Singer (2002) and
by Hartmann (2002), that allow one to reduce the inverse problem for arbitrary
linear algebraic groups over $C(x)$ to finding {\em equivariant} differential equations  with  given connected Galois
groups over an arbitrary finite Galois extension $K$ of $C(x)$.  In Section 3, we develop the
necessary group theory and give  criteria for a differential equation to have a given semisimple
group as its Galois group.  
 In Section 4, we produce equations for the groups described above. In
Section 5, we describe an alternate construction for groups of the form $ W \semi \SL_2$, $W$ a
finite group.

\section{Equivariant Equations}

Mitschi and Singer (2002), and Hartmann (2002), have shown how to reduce the inverse problem for
general groups to the inverse problem for groups of the form $H \ \semi G^0$, where $H$ is a finite
group and $G^0$ is a connected group defined over $C$. We shall assume that we are given a finite
extension $K$ of $C(x)$ with Galois group $H$ and we wish to find a Picard-Vessiot extension $E$ of
$C(x)$ containing $K$ with Galois group $H \ \semi G^0$ such that the Galois action of $H \
\semi G^0$ on $K$ factors through the given action of $H$ on $K$.  To attack this problem, the
authors introduced the notion of an equivariant equation, which we now review.
\\
[0.1in]
Let $K$ be a finite Galois extension of $C(x)$ with Galois group $H$ and let $V$ be a vector space
over $C$ that is also a right $H$-module. Notice that this  gives a right action of $H$ on $V(K)
= K\otimes V$ given by $f\otimes v \mapsto f\otimes v\cdot h$ for $h \in H, f \in K$ and $v \in V$.
We will again denote this action by $w \mapsto w h$ for $h \in H$ and $w \in V(K)$. The group $H$
can be seen to also act on the left on $V(K) $ {\it via} an action defined by $h(f\otimes v) = hf\otimes
v$ for $f
\in K$ and $v \in V$. We say that an element $w \in V(K)$ is {\em equivariant} if $hw = w h$ for all
$h \in H$ (Mitschi and Singer, 2002,  Definition 6.1; Hartmann, 2002, Definition 3.5).
\\
[0.1in]
Let us now consider a semidirect product $G = H \ \semi G^0$ of the finite group $H$ and a connected
linear algebraic group $G^0$ (defined over $C$)  with multiplication given by $(h_1,g_1)(h_2,g_2) =
(h_1h_2, h_2^{-1}g_1h_2g_2)$. Let $\calG$ be the Lie algebra of $G^0$. For any $h \in H$ the map $g
\mapsto h^{-1}gh$ from $G^0$ to $G^0$ can be lifted to a map of $\calG$ to $\calG$ which we shall
also denote by $A \mapsto h^{-1}Ah$. In this way we may consider $\calG$ as a {\em right}
$H$-module. With this convention, we may speak of equivariant elements
 of $\calG(K)$. In concrete terms, an element $A\in\calG (K)$ is equivariant if, for 
 any $h\in H$, the result of applying $h$ (as an element of the Galois group of $K$) to the entries
og $A$ is the same as conjugating $A$ by $h^{-1}$.  We will say that a differential equation $Y' = AY$ is equivariant if $A\in \calG(K)$ is
equivariant. 
We will throughout the paper use the notation $Y'$ for $
\partial Y$, where $
\partial$ is the unique extension of $d/dx$ on $K$. Using this notion of equivariance, we
have the following criterion (Mitschi and Singer, 2002, 
Proposition 6.3; Hartmann, 2002, Proposition 3.10):
\begin{prop} Let $G$ and $K$ be as above and let $A $ be an equivariant element of  $\calG(K)$ such
that the Picard-Vessiot extension $E$ of $K$ corresponding to the equation  $Y' = AY$ has Galois
group $G^0$. Then $E$ is a Picard-Vessiot extension of $C(x)$ with Galois group $G$.
\end{prop}
To apply this result, we will need ways of constructing equivariant elements $A$ of $\calG(K)$ and
criteria to ensure that the equation $Y'=AY$ has the desired Galois group over $K$. 
The remainder of this section is devoted to the first task and the next sections  to the second task.
In the following, we shall think of $C$ as embedded in the complex numbers and denote by $C\{t\}$
the subring of convergent power series of $C[[t]]$,  by $C(\{t\})$  its quotient field, and by $\P1$ the projective line $\P1 (C)$.
\begin{lem} \label{lem1} Let $\pi: \curve \rightarrow  \P1$ be a covering   
of the projective line by a curve $\curve$  with function field  $K$, such that $C(x)\subset K$ is
induced by $\pi$.  There exists a computable set of points $\calS
\subset \P1$ such that the following is true: Given
\begin{enumerate}
\item an integer $M$,
\item points $p_1, \ldots p_r \in \curve$ with $\pi(p_i) \notin \calS$ and $\pi(p_i) \neq \pi(p_j)$
for $i \neq j$,
\item local parameters $t_i$ at $p_i$  and
 \item elements $A_1, \ldots , A_r \in \calG(C)$,
\end{enumerate} there exists an equivariant $A \in \calG(K)$ such that
 at each $p_i$,  we have
\[
A = A_it^M + t^{M+1}(B_i(t))
\]
where $t = t_i$ and $B_i(t) \in \calG(C\{t\}).$
\end{lem}
\begin{proof} We can consider $\calG(C)$ as a {\em left} $H$-module under the action $v \mapsto
hvh^{-1}$ for any $h \in H$. We define an action of $H$ on $\calG(K)$ by the formula $h(f\otimes v)
= h(f)\otimes hvh^{-1}$ for $f\in K$, $v\in\calG (C)$ and $h\in H$. This action satisfies $h(aw) = h(a) h(w)$ for all $h \in H$, $a\in K$,
 $w \in \calG(K)$ and so, by a result of Kolchin and Lang (1993, see Exercises 31 and 32,
p.~550), one can construct an invariant basis of $\calG(K)$ over $K$, that is, a basis $\tildee_1,
\ldots ,\tildee_s$ such that $h(\tildee_i) = \tildee_i$ for $i = 1, \ldots ,s$ for all $h\in H$.  This basis is an
equivariant basis in the above sense, that is, $h\tildee_i = \tildee_i h$ for all $h\in H$. Fix a
basis $e_1, \ldots ,e_s$ of $\calG(C)$ and define $\ B \in \GL_s(K)$  such that $(\tildee_1, \ldots
,\tildee_s) = (e_1, \ldots ,e_s) B$. For
 any $f_1,
\ldots , f_s
\in C(x)$, 
$\sum_{i=1}^sf_i\tildee_i$ is an equivariant element of $\calG(K)$. We shall now show how one can
select the $f_i$ so that the conclusions of the Lemma are satisfied. Let $\calS$ be the image under
$\pi$ of those points $p \in
\curve$ satisfying at least one of the following conditions:
\begin{enumerate}
\item $p$ is a singular point of $\curve$ or is a ramification point of $\pi$, or
\item $p$ is a pole of an entry of  $B$, or
\item $\{\tildee_1(p), \ldots , \tildee_s(p)\}$ fails to be a basis of $\calG(C)$, {\it i.e.}, $\det(B(p))
= 0$. 
\end{enumerate} Note that condition 1.~implies that we may select $t = x-\pi(p)$ to be a local
coordinate for any point $p$ with $\pi(p) \notin \calS$, if $\pi(p)$ is finite, and $t = 1/x$ if
$\pi(p)$ is infinite. We shall use these local coordinates.  From  conditions 1. 2. 3. above, we see that at
each $p_i$ with $\pi(p_i) \notin \calS$, there exist coefficients $c_{i,j}
\in C$ such that $A_i =
\sum_{j=1}^s c_{i,j}\tildee_j(p_i)$.  Let $f_j \in C(x)$ satisfy $f_j = c_{i,j}t^M + t^{M+1}
b_{i,j}$  where $b_{i,j} \in C\{ t\}$ when written in local coordinates $t = t_j$ at the point
$p_j$. We then have that $A =
\sum_{j=1}^s f_j\tildee_j$ satisfies the conclusion of the Lemma.\end{proof}
\begin{cor}\label{cor1}  Let $\pi: \curve \rightarrow  \P1$ be a covering of the projective line by
a curve $\curve$ with function field $K$ such that $C(x)\subset K$ is induced by $\pi$.  There
exists a computable set of points $S
\subset \P1$ such that the following is true: Given
\begin{enumerate}
\item  integers $M<N$,
\item points $p_1, \ldots p_r \in \curve$ with $\pi(p_i) \notin \calS$ and $\pi(p_i) \neq \pi(p_j)$
for $i \neq j$,
\item local parameters $t_i$ at $p_i$  and
 \item for each $i = 1, \ldots ,r$ elements $A_{i,M}, \ldots , A_{i,N} \in
\calG(C)$,
\end{enumerate} there exists an equivariant $A \in \calG(K)$ such that
 at each $p_i$,  we have
\[
A = A_{i,M}t^M + \ldots + A_{i,N}t^N+  t^{N+1}(B_i(t))
\]
where $t = t_i$ and $B_i(t) \in \calG(C\{t\})$
\end{cor}
\begin{proof} Let $\calS$ be as before.  We will proceed by induction on $N-M$. Assume that we have
found an equivariant $A_0 \in \calG(K)$ such that
 at each $p_i$,  we have
\[
A_0 = A_{i,M}t^M + \ldots + A_{i,N-1}t^{N-1}+  t^{N}(B_i(t))
\]
where $t = t_i$ and $B_i(t) \in \calG(C\{t\})$. Let $C_i$ be the coefficient of $t^N$ in the
expansion of $A_0$ at $p_i$. Using Lemma~\ref{lem1}, we can find an equivariant $\tilde{A} \in
\calG(K)$ such that
 \[
\tilde{A} = (A_{i,N} - C_i)t^N+  t^{N+1}(\tilde{B}_i(t))
\]
where $t = t_i$ and $\tilde{B}_i(t) \in \calG(C\{t\})$. The element $A= A_0 + \tilde{A}$ satisfies
the conclusion of the Corollary.
\end{proof}
\begin{ex}{\em The group ${\bf Z}/2{\bf Z} \ \semi \SL_2$}\\
[0.1in] {\em We shall illustrate the above results for this group  where  $ h = -1 \in {\bf Z}/2{\bf
Z}$ acts on $\SL_2$ by sending a matrix to the transpose of its inverse. Note that the action of
this element on $\sl_2$ sends a matrix to the negative of its transpose. Let $K = C(x,\sqrt{x})$
with Galois group $H \simeq {\bf Z}/2{\bf Z}$. The elements
\begin{eqnarray*}
\tildee_1& = &\left(
\begin{array}{cc}
\sqrt{x}&0\\
0&-\sqrt{x} 
\end{array}
\right)\\
\tildee_2 & = & \left(
\begin{array}{cc} 0&\sqrt{x}\\
\sqrt{x}& 0 
\end{array}
\right)\\
\tildee_3& = &\left(
\begin{array}{cc} 0&-1 +\sqrt{x}\\
1+\sqrt{x}& 0
\end{array}
\right)
\end{eqnarray*}
form an equivariant basis of $\sl_2(K)$. We now will construct an equivariant element $A$ of
$\calG(K)$ with the following prescribed principal parts at the points $(4,2), (9,3) \mbox { and }
(16,4)$ of the curve  $y^2 - x = 0$ (we will see in Section 4 that 
the equation $Y' = AY$ is then an equivariant equation with Galois group $\SL_2$ over $K$). 
\begin{eqnarray*}
\mbox{At }  p_0=(4,2), \mbox{ with } t=x-4,\ A &=&
\frac{\left(
\begin{array}{cc} 0&0\\
1&0 
\end{array}\right)}{t^2} +
\frac{\left(
\begin{array}{cc} 0&1\\
0&0 
\end{array}\right)}{t} (\mbox{ terms
 involving }  t^j, \ j \geq 0)\\
\mbox{At }  p_1=(9,3), \mbox{ with } t=x-9,\ A &=&
\frac{\left(
\begin{array}{cc} \sqrt{2}&0\\
0&-\sqrt{2}
\end{array}\right)}{t} + \mbox{ terms involving }  t^j, \ j \geq 0\\
\mbox{At }p_2=(16,4), \mbox{with } t=x-16, A &=&
\frac{\left(
\begin{array}{cc} 0&0\\
1&0 
\end{array}\right)}{t} + \mbox{ terms involving } t^j, \ j \geq 0
\end{eqnarray*}
 A calculation shows that 
the following rational functions $f_i$ yield the desired result for $A:= f_1\tildee_1 +
f_2\tildee_2+f_3\tildee_3$ 

\begin{eqnarray*}
f_1 & = & -{\frac {\sqrt {2}}{105} \frac{\left( x-4 \right)  \left( x-16 \right) }{x- 9}}\\
f_2 & = &-{\frac {1}{240}}\,{\frac { \left( x-9 \right)  \left( x-16 \right) }{
 \left( x-4 \right) ^{2}}} +{\frac {311}{28800}}\,{\frac { \left( x-9
 \right)  \left( x-16 \right) }{x-4} -\frac{3}{672}\,{\frac { \left( x-9 \right)  \left( x-4
 \right) }{x-16}}}\\
f_3 & = & {\frac {1}{120}}\,{\frac { \left( x-9 \right)  \left( x-16 \right) }{
 \left( x-4 \right) ^{2}}}-{\frac {43}{7200}}\,{\frac { \left( x-9
 \right)  \left( x-16 \right) }{x-4}+\frac{1}{168}{\frac { \left( x-9 \right)  \left( x-4
 \right) }{x-16}}.}
\end{eqnarray*}
\hfill \QED }
\end{ex}

\section{Group Theory and its Differential Consequences} Throughout this section  $H$ stands for any
(not necessarily finite) subgroup of a given linear algebraic group $G$. In what follows we shall show that for certain algebraic groups $H$ defined over $C$, we can construct a differential equation, defined over $C(x)$ such that the Galois group of this equation over ${\bf C}(x)$ is $H({\bf C})$. By Theorem 2, section VI.3, of Kolchin (1973) this will imply  that the Galois group of the equation over $C(x)$ is $H(C$).

\smallskip  
\noindent The principal tool that we shall use is the following result:

\begin{lem}\label{lem2} Let $G\subset\GL_n(C)$ be a connected semisimple algebraic group of rank
$\ell$
with Lie algebra $\calG$ and let $\Ad : G
\rightarrow \GL(\calG)$ be the adjoint representation.  Let $H$ be an algebraic subgroup of $G$ and
assume:
\begin{enumerate}
\item $H$  acts reductively on $\calG$ {\it via} the adjoint representation,
\item $H$ contains an element having at least $\ell$ multiplicatively independent eigenvalues,
\item $H$ contains an element $u$ such that ${\rm Ad}(u)$  is a unipotent element with an $\ell$ dimensional
eigenspace corresponding to the eigenvalue 1.
\end{enumerate} Then $H=G$.
\end{lem}

\begin{proof} Let $\calH$ be the Lie algebra of $H$. It is enough to show that $\calH = \calG$.  To
do this we will use the following result (see Exercise 5, p. 246 in Bourbaki, 1990). {\em  If $\calG$ is a
semisimple Lie algebra and $\calH$ is a subalgebra acting reductively on $\calG$ {\it via} the adjoint
representation, of the same rank as $\calG$, and containing a principal $\sl_2$-triple of $\calG$, 
then $\calG = \calH$.}  We shall show that each of the conditions of the Lemma implies the
corresponding condition of this latter result.\\
[0.2in]
If $H$ acts reductively on $\calG$ then so does $\calH$. \\
[0.2in]
Let $h$ be  an element of $H$ having $\ell$ multiplicatively independent eigenvalues.  Then $h^i$ has
the same property for all $i >0$, so we may assume that $h \in H^0$.  We may write $h =h_sh_u$ where
$h_s$ and $h_u$ are the semisimple and unipotent parts of $h$, respectively. Since the eigenvalues
of $h$ and $h_s$ coincide and $h_s \in H^0$, we may assume that  $h$ is semisimple and that $h$ lies
in some maximal torus $T$ of $H$. We may assume that $T$ is a subgroup of diagonal matrices and that
the $\ell$ multiplicatively independent eigenvalues of $h$ are the first $\ell$ entries on the diagonal
of
$h$.  If $\chi_i$ denotes the character that picks out the $i^{th}$ element on the diagonal, we see
that no nontrivial power product of $\chi_1, \ldots ,\chi_\ell$ is trivial on $h$.
 Therefore $\chi_1, \ldots ,\chi_\ell$ are multiplicatively independent on $T$ and so the dimension of
$T$ is greater than or equal to $\ell$, that is, the rank of $\calH$ is $\ell$.\\
[0.2in]
Let $u \in H$ satisfy the property that $\Ad(u)$ is unipotent with an $\ell$-dimensional eigenspace
corresponding to the eigenvalue 1. Since $\Ad(u) = \Ad(u_s)\Ad(u_u)$ where $u_s$ and $u_u$ are
semisimple and unipotent, we have $\Ad(u_s) = 1$ and we can replace $u$ with $u_u$ and assume that $u$
is unipotent. If we let $n =
\log(u)$ we have that $\{\exp(an) \ | \ a \in C\}$ is the unique smallest closed subgroup of $H$
containing $u$ and  that $n \in \calH$ (see Lemma C, p. 96 and Exercise 11, p. 101 of
Humphreys, 1975).  This implies that $\ad(n) = \log(\Ad(u))$ and so a calculation implies that the
dimension of the nulspace of $\ad(n)$ is $\ell$. Since $n$ is nilpotent, the Jacobson-Morozov Theorem
(see Bourbaki, 1990, p. 162) implies that $n$ is contained in an $\sl_2$-triple in $\calH$ and
furthermore this triple is principal in $\SL_n$ ({\it loc. cit.}, p. 166).
\end{proof}
The above lemma gives us the following criterion to ensure that a differential equation has a given
semisimple group as its Galois group.

\begin{prop} \label{prop2} Let $G \subset \GL_n(C)$ be a connected semisimple 
linear algebraic group
of rank $\ell$ with Lie algebra $\calG$ and $\curve$ a curve with function field 
$K \supset C(x)$.
Let $Y' = AY$ be a differential equation with $A \in
\calG(K)$. Let $H\subset\GL_n(C)$ be the Galois group of $Y' = AY$ over $K$ with respect to a 
given fundamental solution $y\in G(K)$, and assume that
\begin{enumerate}
\item $H$ is reductive,
\item There exists a point $p_1 \in C$ such that in terms of some local 
coordinate $t$ at $p_1$, we can expand $Y'=AY$  as 
\[
\frac{dY}{dt} = (\frac{A_1}{t} + B_1(t))Y
\]
where $A_1$ is semisimple and has $\ell$ eigenvalues that are $\bf 
Z$-independent mod $\bf Z$ and
$B_1(t) \in C\{t\}$.
\item There exists a point $p_2 \in C$ such that in terms of some local 
coordinate $t$ at $p_2$, we can expand $Y'=AY$  as
\[
\frac{dY}{dt} = (\frac{A_2}{t} + B_2(t))Y
\]
where $A_2$  is nilpotent and $\ad(A_2)$ has a kernel in $\calG$ of dimension 
$\ell$, and $B_2(t) \in C\{t\}$.
\end{enumerate} Then $H=G$.
\end{prop}
\begin{proof} We know by Proposition 2.1 of Mitschi and Singer (2002) that $H$ 
is an algebraic subgroup of $G$. We
shall show that it  satisfies the hypotheses of Lemma~\ref{lem2}. Clearly 
hypothesis 1.~is
satisfied. \\
[0.2in]
To see that hypothesis 2.~of Lemma~\ref{lem2} is satisfied, note that the 
present hypothesis
2.~implies  that the distinct eigenvalues of $A_1$ do not differ by integers.  
This implies that the
equation $Y'=AY$ is equivalent (over $C((t))$ and even $ C(\{t\}))$ to the 
equation  $Y'
= ({A_1}/{t})Y$ whose monodromy matrix at $p_1$ is  $e^{2\pi i A_1}$ (see 
Section 3 of
Babbitt and Varadarajan, 1983,  or Sections 3.3. and 5.1.1 of
van der Put and Singer, 2003). Note that  $e^{2\pi i A_1}$  has at least $\ell$multiplicatively independent eigenvalues and that this element belongs to $H$.
Since $K$-equivalent differential equations have conjugate Galois groups there 
exists an element $h$
of $H$ which is conjugate (in $\GL_n(C)$) to $e^{2\pi i A_1}$, and hence 
satisfies hypothesis 2.~of Lemma~\ref{lem2}.\\
[0.2in]
To verify hypothesis 3 we must argue in a more careful way. For this we use the 
results of Section 8
of
Babbitt and Varadarajan (1983). Since $\ad(A_2)$ is nilpotent, the spectral 
subspaces ${\bf g}_{\lambda}$ of
$\ad(A_2)$ corresponding to all positive integers $\lambda$ are zero ({\it loc. 
cit.}  Section 8.5).  
Therefore Proposition 8.5 and Theorem 9.5 of Babbitt and Varadarajan (1983) 
imply that there exists 
$g \in G(C\{t\})$ such that the gauge transform of $Y'=AY$ by $Y = gZ$ is 
$Z' =( {A_2}/{t})Z$, and that $A_2$ lies in $\calG$. We note that $Z' =( 
{A_2}/{t})Z$ has a fundamental solution matrix
of the form $z_0 = e^{A_2\log t}$ and that $z_0$ is an element of $G \subset 
\GL_n$ since $A_2\in \calG(C)$.
We therefore have $y=gz_0\gamma$, where $\gamma$ is a constant matrix belonging 
to $G$ since $y$, $g$,  $z_0$ all are elements of $G$. The monodromy matrix at 
$p_2$ with respect to $y$ is  $u=\gamma^{-1} e^{2\pi i A_2} \gamma$. It is an 
element of the local Galois group of $Y'=AY$ over $C(\{t\})$, hence of the 
global Galois group $H$ over $K$ with respect to $y$. Since $\gamma\in G(C)$, 
the matrices $\Ad(u)$ and $\Ad(e^{2\pi i A_2})$ are conjugate in $\GL(\calG)$, 
which implies that $Ad(u)$ has the desired property.  \end{proof} 

%
\noindent From Section 2 it is clear that we should have no trouble fulfilling hypotheses 2.~and 3.~of
Proposition~\ref{prop2}. The difficulty arises in trying to ensure that hypothesis~1.~is
satisfied.
We shall use local properties of the differential equation. As a first step we
will derive a condition on the behavior of a differential equation at {\em one} singular point to
ensure that the Galois group is reductive. We will see that this  will only work for
$\SL_n$ and $\Sp_n$.  We will then give a more general criterion that involves the local behavior
at several points and this will apply to a larger class of groups.\\[0.1in]
Before we describe criteria in terms of local properties of the differential equation 
that ensure that  hypothesis 1. of Proposition~\ref{prop2} is satisfied, we will recall the facts we need relating 
the local Galois groups and the global Galois group. \\
[0.2in]
 Let $\pi: \curve \rightarrow \P1$ be a nonsingular curve over the projective line and  $C(x)
\subset K$ the corresponding inclusion of function fields.
 Let $p \in \curve$ and assume that $\curve$ is not ramified at $p$ and that $\pi(p) \neq \infty$.
What follows can be developed without these assumptions but they simplify the exposition and will
hold in our applications. We can  embed $K$ into $C((t)),
\ t = x-\pi(p)$, by expanding each element of $K$ as a series in $t$. We shall identify $K$ with its
image and write $K \subset C((t))$.  In fact, we have that $K \subset C(\{t\}) \subset C((t))$. Any
differential equation $Y' = AY, \ A \in
\gl_n(K)$ can be considered as a differential equation over $C((t))$ and so we can form a
Picard-Vessiot extension $E$ of $C((t))$ corresponding to this equation.   Let $y$ be a fundamental
solution  of $Y' = AY$ having entries in $E$, and let $K(y)$ and $C(\{t\})(y)$ denote the
fields generated by the entries of $y$ over $K$ and $C(\{t\})$ respectively. We see that $K(y)$ and $C(\{t\})(y)$ are
Picard-Vessiot extensions for $Y' = AY$ over $K$ and $C(\{t\})$ respectively. We denote by $G$, $G_{conv} \mbox{ and } G_{form}$ the Galois groups
of $K(y)$ over $K$, of $C(\{t\})(y)$ over $C(\{t\})$ and of $E$ over $C((t))$ respectively.  One 
easily checks that there are natural injections $G_{form}
\hookrightarrow G_{conv}
\hookrightarrow G$ and that the actions of the former two groups on the solution space of the
differential equation coincide with their actions as embedded subgroups of $G$.  These
considerations lead to:
\begin{lem}\label{lem3} Let $\curve$ be  a curve with function field $K \supset C(x)$ and let $Y' =
AY$ be a differential equation with coefficients in $K$. If there exists a point $p \in \curve$ as above such that the equation  is irreducible
over $C(\{t\})$, then it is irreducible over $K$ and its Galois group $G$ is reductive. In
particular, if it is irreducible over $C((t))$, then $G$ is
reductive.
\end{lem}
\begin{proof} A differential equation is irreducible over a differential field with algebraically
closed field of constants if and only if its Galois group  acts irreducibly on the solution space of
the equation in a Picard-Vessiot extension. If $Y' =AY$ is irreducible over $C(\{t\})$ then
$G_{conv}$ acts irreducibly on the solutions space. Since $G_{conv} \hookrightarrow G$, we have that
$G$ acts irreducibly on this space and so the equation is irreducible over $K$.  We can furthermore
conclude that $G$ is reductive since it has an irreducible faithful representation. The final
statement follows in a similar manner.
\end{proof}

\noindent It is much easier to show that a differential equation  is irreducible over $C((t))$ than to show it
is irreducible over $C(\{t\})$.  Regrettably, from our point of view, irreducibility over
${C((t))}$ puts severe restrictions on the Galois group  of $Y' = AY$ over
${C(x)}$. One can deduce from  Remark 3.34 of van der Put and Singer (2003) that if $Y'=AY$ is
irreducible over $C((t))$, then
\begin{enumerate}
\item  the identity component $G_{form}^0$ of $G_{form}$ is a torus,
\item as a $G_{form}^0$-module, the solution space is the sum of one dimensional invariant subspaces
corresponding to distinct characters of $G_{form}^0$, and
\item there is an element $\gamma \in G_{form}$ whose action on $G_{form}^0$ by conjugation
cyclically permutes the characters of $G_{form}^0$.
\end{enumerate} Katz (1987, 3.2.8 and 3.2.9)  has shown that a connected
algebraic subgroup  of $\SL_n$, containing a closed subgroup  satisfying the
properties 1., 2., and 3. of $G_{form}$ above must be of the form $\prod G_i$ where each $G_i$ is either
$\SL_{n_i}$ or $\Sp_{n_i}, n_i$ even in the latter case,  and the $n_i$ are pairwise relatively
prime. Katz further shows that the $n$-space (in our case the solution space) can be written as a
tensor product $\otimes V_i$ of representations of these groups where each $V_i$ is the standard or
contragredient representation of $G_i$ if $G_i = \SL_{n_i}$ or the standard representation of $G_i$
if $G_i = \Sp_{n_i}$.\\
[.2in]
 Nonetheless, Lemma~\ref{lem3} together with Proposition~\ref{prop2} will allow us to construct 
equations $Y' = AY$ 
 having  Galois group $\SL_n$ or $\Sp_{2n}$. These two results yield the following criteria:

\begin{prop} \label{prop3} Let $G \subset \SL_n(C)$ be a connected simple linear algebraic group of
rank $\ell$ with Lie algebra $\calG$ and $\curve$ a curve with function field $K \supset C(x)$. Let
$Y' = AY$ be a differential equation with $A \in
\calG(K)$. Let $H\subset\GL_n(C)$ be the Galois group of $Y' = AY$ over $K$ with respect to a given fundamental solution $y\in G(K)$ and assume that
\begin{enumerate}
\item There exists a point $p_0 \in \curve$ such that the equation $Y' = AY$ has a unique slope of the
form $\frac{a}{n}, \
(a,n) = 1$.
\item There exists a point $p_1 \in \curve$ such that in terms of some local coordinate $t$ at
$p_1$, we can expand $Y'=AY$ as
\[
\frac{dY}{dt} = (\frac{A_1}{t} + B_1(t))Y
\]
where $A_1$ is semisimple and has $\ell$ eigenvalues that are $\bf Z$-independent mod $\bf Z$ and
$B_1(t) \in \gl(C\{t\})$.
\item There exists a point $p_2 \in \curve$ such that in terms of some local coordinate $t$ at
$p_2$, we we can expand $Y'=AY$ as
\[
\frac{dY}{dt} = (\frac{A_2}{t} + B_2(t))Y
\]
where $A_2$ is nilpotent and the kernel of $ad(A_2)$ in $\calG$ has dimension $\ell$, and $B_2(t) \in \gl_n(C\{t\})$.
\end{enumerate} Then $H=G$. Furthermore, if this is the case, then $G$ must be either $\SL_n$ or
$\Sp_{2n}$.\end{prop}
\begin{proof} We refer to Babbitt and Varadarajan (1983), Katz (1987) or
van der Put and Singer (2003) for the definition and properties of the slopes of a differential equation at a
singular point.  From (2.2.8) of Katz (1987) or Remark 3.34 of van der Put and Singer (2003), one
sees that hypothesis 1.~above implies hypothesis 1.~of Proposition \ref{prop2}. The last statement
follows from the discussion preceding the statement of this proposition.
\end{proof}

\noindent To give irreducibility criteria that apply to other groups we shall show how one can
compare
 the local behavior at several points to ensure irreducibility. These criteria will assume that at
several points the 
identity component of the local formal Galois groups is a maximal torus of the global Galois group (that is, the Galois group over the function field $K$ of $\curve$) and they will give conditions on
elements normalizing these tori to ensure irreducibility.  We therefore start with
the following definition.

\noindent We continue the convention that 
$\pi: \curve \rightarrow \P1$ is a nonsingular curve over the projective line,
 $C(x) \subset K$ is the corresponding inclusion of function fields and 
  $p \in \curve$ with $\curve$  not ramified at $p$ and that $\pi(p) \neq \infty$. 

\begin{defin}\label{toricdef} Let $G$ be a connected linear algebraic group with Lie algebra $\cal G$.  
Let $Y' = AY$ be a differential equation with $A\in \calG(K)$. We say that $p \in \curve$ is a 
{\em maximally toric point} for $Y' = AY$ (with respect to $G$) if the connected
component 
$G_{form}^0$ of 
the local formal Galois group $G_{form}$ at $p$ is a maximal torus in $G$.
\end{defin}

\noindent The work of Katz quoted above implies that if $p \in \curve$ is a point such that the
equation has a unique slope of the form 
${a}/{n}, (a,n) = 1$, then $p$ is a maximally toric point but we shall see that not all
maximally toric points need arise in this way.\\[0.1in]
\noindent Let $Y' = AY$ be a differential equation  as in Definition~\ref{toricdef} and let 
$p$ be a maximally toric point for this equation and let $G_{form}$ be the
formal local Galois groups at $p$. Theorem 11.2 of van der Put and Singer (2003) implies 
that $G_{form}/G_{form}^0$ is a finite cyclic group. Since $G_{form}^0$ is a maximal torus of
$G$, we may identify the generator $g$ of  $G_{form}/G_{form}^0$ with an element
of the
 Weyl group of $G$. Regrettably, we do not see how to do this in a canonical way so that we can
compare the images of elements $g$ for different maximally toric points $p$ of $\curve$. Nonetheless,  assume that a (and therefore any) maximal torus of $G$ has $m$ weight
spaces in the representation of $G$ on the solution space of $Y' = AY$. 
Since the element $g$  permutes the  weight spaces, it can be considered as an
element of $\gS_m$, the permutation group on $m$ elements (again in a noncanonical way). 
The key fact is that although the image of $g$ in $\gS_m$ may not be uniquely defined, all such
images are conjugate in $\gS_m$ since  $g$ is determined up to conjugation in $G$.  We refer to
this $\gS_m$ conjugacy class as the {\em permutation conjugacy class} at the toric point $p$.
We now make the following definition.

\begin{defin} Let $\gS_m$ be the permutation group on $m$ elements, and 
  let $\{C_1, \ldots ,C_t\}$ be a collection of conjugacy classes in  $\gS_m$. 
   We say that the set $\{C_1, \ldots ,C_t\}$ is {\em strictly transitive}
 if for any choice $\tau_i \in C_i, i = 1, \ldots t$, the subgroup of $\gS_m$ generated by 
 $\{\tau_1, \ldots ,\tau_t\}$
 acts transitively.
\end{defin}
Since the conjugacy class of an element of $\gS _m$ is determined by the type of the partition
on $\{1,\ldots,m\}$ given by its cycle structure, one can see that 

\begin{lem}\label{transitive} The set $\{C_1, \ldots ,C_t\}$ of conjugacy classes in $\gS_m$  is 
strictly
transitive if and only if the  following holds for some (and therefore any) set of representatives
$\{\sigma_1,
\ldots ,
\sigma_t\}$ with $\sigma_i \in C_i$:
 for any $i, 1\leq i \leq m-1$, there is an element $\sigma_j$ leaving no set of
cardinality $i$ invariant.
\end{lem}

\noindent For example, for $m=6$ the singleton set $\{\overline{(123456)}\}$ and the set 
$\{\overline{(123)(456)}, \overline{(1234)(56)}\}$ are
strictly transitive sets of conjugacy classes (where $\overline{\sigma}$ denotes the conjugacy
class of $\sigma$).  The set $\{\overline{(123)(456)}, \overline{(1)(45)(236)}\}$ is
not strictly transitive since each permutation leaves a set of  $3$ invariant. We
are now ready to state the following criterion, which generalizes Proposition~\ref{prop3}
\begin{prop}\label{prop4} Let $G \subset \GL_n(C)$ be a connected simple linear algebraic group of
rank $\ell$ with Lie algebra $\calG$ and $\curve$ a curve with function field $K \supset C(x)$. Let
$Y' = AY$ be a differential equation with $A \in
\calG(K)$. Let $H$ be the Galois group of $Y' = AY$ over $K$ and assume that
\begin{enumerate}
\item There exist  maximally toric points $p_1, \ldots , p_t \in \curve$ for the equation $Y' = AY$
such that the corresponding conjugacy classes form a strictly transitive set.
\item There exists a point $p_1 \in \curve$ such that in terms of some local coordinate $t$ at
$p_1$, we can expand $Y'=AY$ as
\[
\frac{dY}{dt} = (\frac{A_1}{t} + B_1(t))Y
\]
where $A_1$ is semisimple and has $\ell$ eigenvalues that are $\bf Z$-independent mod. $\bf Z$ and
$B_1(t) \in \gl_n(C\{t\})$.
\item There exists a point $p_2 \in \curve$ such that in terms of some local coordinate $t$ at
$p_2$, we can expand $Y'=AY$ as
\[
\frac{dY}{dt} = (\frac{A_2}{t} + B_2(t))Y
\]
where $A_2$ is nilpotent and the kernel of $\ad(A_2)$ has dimension $\ell$ and $B_2(t) \in \gl_n(C\{t\})$.
\end{enumerate} Then $H=G$. 
\end{prop}
\begin{proof} As in the proof of Proposition~\ref{prop3}, we need only to show that the first
condition guarantees that the Galois group $H$ acts irreducibly on the solution space $V$ and so that  $H$ is
reductive. Let $T$ be a maximal torus of $H$ and assume that $T$ has $m$ distinct weight
spaces in the solution space of $Y' = AY$.  If $W$ is a proper, nontrivial
 $H$-invariant subspace of  $V$ then we can write $W$ as a sum 
of $i, 1\leq i \leq m-1,$ weight spaces.  Since each local formal Galois
group $G_{form,j}$ at $p_j$ leaves $W$ invariant, we can conclude that  for each $j$, any generator
$\sigma_j$ of the group $G_{form,j}/ G_{form,j}^0$ leaves a set of $i$ weight spaces stable.
Therefore, Lemma~\ref{transitive} implies that the set of associated conjugacy classes is not strictly
transitive, a contradiction.  Therefore $V$ is an irreducible $H$-module.
\end{proof}

\noindent We note here that if we can satisfy the first condition of Proposition~\ref{prop4}, then the
representation of $G \subset \GL_n$ is severely restricted.   In
particular the Weyl group of $G$
will act transitively on the weights and so the representation will be a minuscule representation (see
Ch.VIII, \S7, No. 3 of Bourbaki, 1990, and Katz, 1987,  p. 48). This means that if $G$ is a simple group it
must be  of type $\Al, \Bl, \Cl, \Dl, E_6, E_7$ and the representations must have highest weight given in the
following list:
\begin{eqnarray*}
\Al (\ell \geq1): & & \obar_1, \obar_2, \dots, \obar_\ell\\
\Bl (\ell \geq 2): & &\obar_\ell\\
\Cl (\ell \geq 2): & & \obar_1\\
\Dl (\ell \geq 3): & & \obar_1, \obar_{\ell-1}, \obar_\ell\\
\esi: & & \obar_1, \obar_6\\
\ese: & & \obar_7\\
{\rm E}_8, {\rm F}_4, {\rm G}_2: & & \mbox{no minuscule weights}
\end{eqnarray*}

\noindent In fact we shall see that there are groups of type $\Bl$ which have no representation with a
strictly transitive set of conjugacy classes in the Weyl group but we are able to apply the
proposition to the rest of the possible types.\\[0.1in]
To apply Proposition~\ref{prop4}, we need to ensure that, given a  group $G$ with Lie algebra $\calG$, we are able to
construct an $A \in \calG(K)$ having prescribed formal Galois groups at given points.  The inverse
problem over $C((t)), t' = 1$, has been solved  and it is known which groups appear as formal Galois groups and 
how one can even effectively construct equations $Y' = AY$, ${A \in \gl_n(C((t)))}$ with allowable formal Galois groups (see
Chapter 11.2 of van der Put and Singer, 2003). We
will
need a small modification of this result to ensure that ${A \in \calG(C((t)))}$ where $\calG$ is the
Lie algebra of a given group.  This is done in Proposition~\ref{local}. We begin with two lemmas.
Lemma~\ref{toric90} is a slight generalization of Hilbert's Theorem 90.
\begin{lem}\label{toric90} Let $E$ be a finite Galois extension of a field $F$ and $G$ a connected linear algebriac group
defined over an algebraically closed field $C \subset F$.  Assume that the Galois group $H$ of $E$ over $F$ is a
finite cyclic group and let $\rho: H \rightarrow G(C)$ be a homomorphism.  
Then there is an element $\eta \in G(F)$ such that $\eta^{\tau} = \eta\cdot \rho(\tau)$ for all $\tau\in H$, where $\eta^\tau$
denotes  the result of applying $\tau$ to the entries of $\eta$.
\end{lem}
\begin{proof}
Since $\rho(H)$ is cyclic, it is contained in a torus $T$ of $G(C)$.  There is a
$C$-isomorphism of $T$ with some power  of the multiplicative group $G_m(C)=C^*$ of $C$ so we can assume that 
$\rho(H) \subset( G_m(C))^r$. The map $\rho$ as a $1$-cocycle defines a cohomology class  $\overline{\rho}\in 
H^1(H,( G_m(C))^r) $ = $(H^1(H, G_m(C)))^r$.
Hilbert's Theorem 90 (see Chapter VI, \S10 of Lang, 1993) implies that $ H^1(H, G_m(C)) = 1$ so $\overline{\rho}$ is  trivial, that is,  $\rho$ is a coboundary. This proves the lemma.
\end{proof} 
The following lemma is a distillation  and modification of the proof of Theorem 11.2 of
van der Put and Singer (2003), included for the convenience of the reader.
\begin{lem}\label{toricautom} Let $T$ be the group of diagonal elements in $\GL_n(C)$ with Lie
algebra $\calT \subset \gl_n(C)$ and $\phi:T\rightarrow T$ an automorphism of order m. There
exists an $A \in \calT(C[t^{-\frac{1}{m}}])$ such that the Galois group of $Y' = AY$ over
$C((t^{\frac{1}{m}}))$ is $T(C)$ and such that $A^{\gamma} = d\phi(A)$ where $A^{\gamma}$ denotes
the matrix resulting from applying the automorphism  of $C((t^{\frac{1}{m}}))$ defined by
$\gamma:t^{\frac{1}{m}} \mapsto e^{\frac{2\pi i}{m}}t^{\frac{1}{m}}$ to the entries of $A$.
\end{lem}
\begin{proof} The map $d\phi:\calT(C) \rightarrow \calT(C)$ is an automorphism of order $m$ and so
has eigenvalues that are roots of unity. Let $W_q \subset \calT(C)$ be a nonzero eigenspace
corresponding to a root of unity $e^{\frac{2\pi i q}{m}}$, where $q$ is an integer $0\leq q <m$ and
let ${\bf b}_{j,q}, 1 \leq j \leq r_q$ be a basis of $W_q$. Defining 
\[ A_q = z^{-\frac{q}{m}}\sum_{j=1}^{r_q} z^{-j-1}{\bf b}_{j,q}\]
one sees that $A_q^\gamma = d\phi(A_q)$.  Let $A = \sum A_q$ where the sum is over all $q$ with
$W_q \neq (0)$.  We claim that $Y' = AY$ satisfies the conclusion of the Lemma. The behavior with 
respect to $\gamma$ follows from the construction of $A$. Since $A
\in
\calT(C((t^{\frac{1}{m}})))$, the Galois group of $Y'=AY$ over $C((t^{\frac{1}{m}}))$ is a subgroup of $T(C)$.
A full solution space for this equation is spanned by \[\{\frac{1}{-\frac{q}{m}
- j}e^{z^{-\frac{q}{m} - j}} {\bf
b}_{j,q}\}
\]
where
$q$ runs over those integers with $W_q \neq (0)$ and $1\leq  j \leq r_q$. To show that the Galois
group is all of $T(C)$ it suffices to show that the elements  $\{e^{z^{-\frac{q}{m} - j}}\}$ form an
algebraically independent set.  The Kolchin-Ostrowski Theorem (Kolchin, 1968) implies
that this is the case if there is no nontrivial relation of the form $\sum a_{j,q} z^{-\frac{q}{m} -
j-1} = \frac{f'}{f}$ where the $a_{j,q}$ are rational numbers and $f \in C((t^{\frac{1}{m}}))$ (see also  Exercise 4, Chapter VI.5 of Kolchin, 1973). Since the
order of $\frac{f'}{f}$ is $ \geq -1$ and the order of a nonzero $\sum a_{j,q} z^{-\frac{q}{m} -
j-1} $ is  $<-1$, we see that no such nontrivial relation can exist. 
\end{proof}
\begin{prop}\label{local} Let $G \subset \GL_n$ be a connected linear algebraic group defined over
$C$ and
$\Gbar \subset G$ a subgroup with $\Gbar^0$ a torus and $\Gbar/\Gbar^0$ cyclic. Let $\calGbar$ be
the Lie algebra of $\Gbar$. There exists an $\Abar \in \calGbar(C[t,t^{-1}])$ such that the Galois group of
$Y' = AY$ over $C((t))$ is $\Gbar$. Furthermore, this $\Abar = \sum_{i=a}^{b} A_it^i$ can be chosen so
that if $A$ is any element of $\calG(C((t))$ such that $A-\Abar = \sum_{i=b+1}^\infty B_it^i$,   then the Galois
group of $Y' = A Y$ over $C((t))$ is also $\Gbar$.
\end{prop}
\begin{proof} 
 We begin by noting that under the assumptions of Proposition \ref{local} there is an element $g \in \Gbar$ of
finite order whose image generates 
$\Gbar/\Gbar^0$ (see  Theorem 8.10 of van der Put and Singer, 1997). To see this let $\gbar$ be an element whose image
generates
$\Gbar/\Gbar^0$. If $m = |\Gbar/\Gbar^0|$, then $\gbar^m \in \Gbar^0$ and so  $\gbar$ is semisimple.
The Zariski closure
$Z$ of the group generated by $\gbar$ has only semisimple elements and is therefore diagonalizable.
Following Theorem 16.2 of Humphreys (1975) we can write $Z$  as the direct product of a torus and a finite group $H$.  Since $Z/Z^0 \rightarrow \Gbar/\Gbar^0$ is surjective, there is
some element of $H$ that maps to a generator of $\Gbar/\Gbar^0$. Let $g$ be this element and
assume $g$ has order $m'$. Note that $m|m'$. We are now going to find elements $\eta, s \in G$ 
such that the element $(\eta s)'\cdot (\eta s)^{-1}$ of $ \calG(k_0)$, with $k_0 = C((t))$, satisfies all but the
last sentence of the above Proposition.\\ [0.1in]
We first select $\eta$. Let $k = C((t^{\frac{1}{m}}))$ and  $k' = C((t^\frac{1}{m'}))$.  We may identify
$k$ with a subfield
of $k'$.  Let $H'$ denote the Galois group of $k'$ over $k_0$.  Since this is a cyclic group
of order $m'$, there is an isomorphism $\rho: H' \rightarrow G(C)$ mapping a generator $\gamma '$ of $H'$ to
$g$. Let $\eta \in G(k')$ be an element, guaranteed to exist by Lemma~\ref{toric90}, that satisfies
$\eta ^{\gamma'} = \eta \cdot g$.\\[0.1in]
We are now ready to select $s$.  Since $\Gbar^0$ is a torus we may identify it with the group $T$ of diagonal elements in some 
$\GL_r(C)$. 
Conjugation by $g$ induces an automorphism of $\Gbar^0$ of order $m$.  Lemma~\ref{toricautom} implies that there is an $\tilde{A} \in
\overline{\calG}(k) \subset \calG(k)$ such that the Galois 
group of $Y' = \tilde{A} Y$ over $k$ is $\Gbar^0$ and $\tilde{A}^\gamma = g^{-1} \tilde{A} g$, where $\gamma$ as before is the $k_0$-automorphism  $t^{\frac{1}{m}} \mapsto e^{\frac{2\pi i}{m}}t^{\frac{1}{m}}$ of $k$. 
We shall now consider $\tilde{A}$ as an element of $\calG(k')$ (since $k\subset k'$) and identify  $\gamma$ with an
element of $H'$, the Galois group of $k'$ over $k_0$. Let $s$ be a $L$-point of  ${\Gbar}^0$ (in a suitable differential extension $L$  of $k$)    satisfying the differential equation $s' = \tilde{A}s$ and such
that $k'(s)$ is a Picard-Vessiot extension of $k'$ with Galois group $\Gbar^0$ for this equation.\\[0.1in]
We now show that $k_0(\eta s)$ is a Picard-Vessiot extension of $k_0$ with Galois group $\Gbar$. As in the proof of Proposition
6.3 of Mitschi and Singer (2002) one sees that 
\begin{enumerate}
\item  The element $(\eta, \eta s)$ (denoted there by $(w, wg)$) satisfies
the equation $Y' = B Y$ where 
\[ B = \left(\begin{array}{cc} \eta' \eta ^{-1} & 0 \\ 0 &  \eta'\eta^{-1} + \eta
\tilde{A}\eta^{-1}\end{array}\right) \ \ .\]
\item Both $\eta'\eta^{-1}$ and  $\eta'\eta^{-1} + \eta\tilde{A}\eta^{-1}$ are in $\Gbar(k_0)$.
\item The Picard-Vessiot extension $k_0((\eta,\eta s))$ of $k_0$ has Galois group $<g>\semi
T(C)$, where $<g>$ is the cyclic group generated by $g$.
\end{enumerate}
We claim that the Picard-Vessiot extension $k_0(\eta s)$ has Galois group $\Gbar$. Note that since
$k_0(\eta s) \subset k_0((\eta ,\eta s))$, the Galois group of $k_0(\eta s)$ is isomorphic to the quotient of 
 $<g>\semi T(C)$ by the subgroup of its elements that leave $\eta s$ fixed.  An element $(a,b)$ of 
$<g>\semi T(C)$ maps $\eta s$ to $\eta sab$.Therefore $(a,b)$
leaves $\eta s$ fixed if and only if $ab = 1$.  The set of such elements is the same as the kernel of
the homomorphism $<g>\semi T(C) \rightarrow <g> T(C) = \Gbar$ that sends $(a,b)$ to $ab$.
Therefore the Galois group of $k_0(\eta s)$ is $\Gbar$.\\[0.1in]
As noted above,  $(\eta s)' (\eta s)^{-1}$ is an element of $\calG(k_0)$. If we
write it as $\sum_{i = p}^{\infty}A_it^i$, then the results of Sections 6 and 7 of
Babbitt and Varadarajan (1983) imply that one can effectively find an integer $q$ such that the canonical
form of the equation $Y' = (\sum_{i = p}^{\infty}A_it^i)Y$ is determined by $\sum_{i = p}^{q}A_it^i$ and so is
its formal Galois group.  Therefore, for $\Abar = \sum_{i = p}^{q}A_it^i$, the Proposition is
proved. \end{proof}
\begin{ex}{\em  We shall illustrate the above proposition for $G = \SL_2$.  Let 
\[\Gbar = \{\left( \begin{array}{cc} a & 0 \\ 0 & a^{-1}\end{array}\right), \ 
\left( \begin{array}{cc} 0 & a \\ -a^{-1} & 0\end{array}\right) \ | \ a \neq 0 \} \ . \] 
The identity component $\Gbar^0$ is a maximal torus and $\Gbar/\Gbar^0$ is the cyclic group of order two. 
The element
\[g = \left( \begin{array}{cc} 0 & 1 \\ -1 & 0\end{array}\right) \]
is an element of order $4$ whose image generates the group $\Gbar/\Gbar^0$. Let $k' = C((t^{\frac{1}{4}}))$ and
\[ \eta = \frac{1}{2}\left( \begin{array}{cc} t^{\frac{1}{4}} + t^{-\frac{1}{4}} & \sqrt{-1}(-t^{\frac{1}{4}} +t^{-\frac{1}{4}}) \\ 
\sqrt{-1}(t^{\frac{1}{4}} - t^{-\frac{1}{4}}) & t^{\frac{1}{4}}+ t^{-\frac{1}{4}}\end{array}\right) \ . \]
If $\gamma'$ is the generator of the Galois group of $k'$ over $C((t))$, then one sees that $\eta^{\gamma'} = \eta g$.  Now define
\[s = \left( \begin{array}{cc} e^{t^{-1/2}} & 0 \\ 0 & e^{-t^{-1/2}}\end{array}\right) \ .\]
One sees that $s$ satisfies the differential equation $s' = \tilde{A}s$ where 
\[\tilde{A} =  \left( \begin {array}{cc}\frac{ -1}{2{t}^{3/2}}&0\\\noalign{\medskip}0&\frac{1}{2
{t}^{3/2}}\end {array} \right)
 \ ,\]
and that $\eta s$ satisfies the  equation  $(\eta s)' = \overline{A}(\eta s)$ where 
\[\overline{A} =  \left( \begin {array}{cc} -{\frac {t+1}{4{t}^{2}}}&{\frac {-
\sqrt{-1} \left( 2\,t-1 \right) }{4{t}^{2}}}\\\noalign{\medskip}\frac{\sqrt{-1}}
{4{t}^{2}}&{\frac {t+1}{4{t}^{2}}}\end {array} \right)  \ .
 \]
The proposition assures us that the equation $Y' = \overline{A}Y$ has Galois group  $\Gbar$. \hfill \square
}
\end{ex}

\section{Equivariant Equations for $\boldsymbol{G^0 = G_1 \cdot \ldots \cdot G_r}$   where each $\boldsymbol{G_i}$ is a 
Simple Group of Type $\boldsymbol{\Al}$, $\boldsymbol{\Cl}$, $\boldsymbol{\Dl}$,
 $\boldsymbol {\esi}$ or
$\boldsymbol{\ese}$} In the previous section, we gave criteria which guarantee that a given equation has
given Galois group.  Here we shall show how one can construct equations meeting
these criteria.  In
this section $G$ is a given linear algebraic group with identity component $G^0$ and $H$, as in
section 2,  denotes a finite subgroup of $G$. 
Let $K$ be a Galois extension of $C(x)$ with Galois group $H$.  We will begin by showing that it
suffices to accomplish this task for {\em simply connected} groups $G_i$.  We shall then show how to
find equivariant equations for simply connected groups of each of the types mentioned above and then
for products of these groups. \\[0.1in]
To reduce the inverse problem in our situation to an inverse problem for 
simply connected groups we proceed as follows.  Let $G = H \semi G^0$ with $G^0$ a semisimple group.  From Theorem 5.1
of Hochschild (1976) there exist a simply connected group $\tilde{G^0}$ and a morphism $\rho: \tilde{G^0} \rightarrow
G^0$ with finite kernel. Theorem 5.5 of 
the same reference implies that every morphism $\sigma :G^0 \rightarrow G^0$ lifts to a {\em
unique} morphism
$\tilde{\sigma}:
\tilde{G^0}
\rightarrow \tilde{G^0}$ such that 
$\rho \tilde{\sigma}= \sigma\rho$. In particular this implies that the action of $H$ on $G^0$ lifts to an action of $H$
on $\tilde{G^0}$ such that  ${\rm {id}} \times \rho: H\semi \tilde{G^0} \rightarrow H\semi G^0$ is a morphism.
Therefore if we can find an $H$-equivariant $ \tilde{A} \in \tilde{\calG}(K)$ (where $\tilde{\calG}$ is the Lie algebra of $\tilde{G^0}$) such that the  Galois group of $Y' = AY$ 
is $ H\semi \tilde{G^0}$, then 
  $d\rho(A)$ is $H$-equivariant and by Proposition 5.3 of Mitschi and Singer (2002) the Galois group of $Y' = d\rho(A)Y$ is $G$.
\subsection{Equivariant Equations  for Simply Connected Groups of Type  $\boldsymbol{\Al}$, $\boldsymbol{\Cl}$, $\boldsymbol{\Dl}$,
 $\boldsymbol {\esi}$ or
$\boldsymbol{\ese}$} In this section we shall apply Propositions~\ref{prop3} and \ref{prop4} to construct equivariant
equations with groups of the above types. In fact we shall show: 
\begin{prop}\label{propsimply} Let $G$ be a simply connected group of type $\Al$,  $\Cl$,  $\Dl$, $\esi$ or
$\ese$. There is a representation
$\rho:G \rightarrow \GL_N$ with associated representation  $d\rho:\calG \rightarrow \gl_N$  of its Lie algebra, 
and elements $B_1, B_2, ... ,B_m \in d\rho(\calG)[t^{-1}, t]$ such that
\begin{enumerate}
\item  for any covering $\pi:\curve \rightarrow \P1$ 
of the projective line by a nonsingular curve $\curve$ with $C(x) \subset K$ the corresponding extension  of  function fields,  
and points $q_1, \ldots q_m \in \curve$ that are not ramification points and such that $\pi(q_i)\ne \infty$ for  $i=1,\ldots,m$,  and
\item any element $B \in d\rho(\calG)(K)$ such that $B = B_i + $(terms of order higher than ${\rm{deg}}(B_i)$) in a local coordinate $t$ at $q_i$, for $i=1,\ldots,m$,
\end{enumerate}
the equation $Y' = BY$ has Galois group $G$ over $K$. Furthermore, at least one of the $q_i$ is an irregular singular point of $Y' = BY.$
\end{prop}
In fact, for  simple groups under consideration,  $m$ can be chosen to be $3$ or $4$.\\[0.1in]
\noindent We  note that once the proposition has been established, a simple application of Corollary~\ref{cor1} yields an equivariant equation.
We will prove Proposition \ref{propsimply} by showing that
one can select the $B_i$ so that the conditions of Propositions~\ref{prop3} or \ref{prop4} are satisfied.  
We first show how conditions 2~and 3~of these propositions can be fulfilled and then turn to condition 1.\\
[0.1in]
{\bf Condition 2.} If $G$ is a linear algebraic group of rank $\ell$, then, by definition,
$G$ contains a torus $T$ of dimension $\ell$. Any torus contains elements $g$ that  generate  each a
Zariski
dense subgroup of $T$. Such an element must be semisimple and have $\ell$ multiplicatively
independent eigenvalues. We may write $g = e^{A_1}$ for some $A_1 \in \calT(C) \subset \calG(C)$. One sees that
$A_1$ will also be semisimple and have $\ell$ eigenvalues that are ${\bf Z}$-independent {\it mod.} ${\bf
Z}$. For example, let $r_1,
\ldots , r_{n-1} \in C$ be ${\bf Z}$-linearly independent {\it mod.} ${\bf Z}$ and let $r_n = -\sum r_i$. 
If $G^0 = \SL_n$ let $A_1 = {\rm diag}(r_1, \ldots, r_n)\in \sl_n(C)$. If $G^0 = \Sp_{2n}$, let $A_1 = {\rm diag}({\rm diag}(r_1, \ldots,
r_n),-{\rm diag}(r_1, \ldots, r_n))\in\sp_n(C)$.    We let $B_1 = {A_1}/{t}$. \\[0.1in]
{\bf Condition 3.} The element $A_2$ will be a principal nilpotent element of the Lie algebra. By Proposition 8, p. 166, of Bourbaki (1990), any 
semisimple Lie algebra of rank $\ell$ contains principal nilpotent elements $u$. These can be constructed by decomposing the algebra as
the sum of a Cartan subalgebra and nonzero root spaces $\g_{\alpha}$ and letting $u = \sum_{\alpha \in \Phi^+}
v_{\alpha}$ where   $v_{\alpha}$ is a nonzero element of  $\g_{\alpha}$ for each positive root $\alpha\in \Phi^+$ ({\it loc. cit.} Proposition 10, p. 168).  For example, in $\sl_n(C)$ we can take the matrix
$u=(a_{i,j})$ where $a_{i,j} = 0$ if $i \geq j$ and $a_{i,j} = 1$  if $i < j$. We let $A_2$ be such
an element and $B_2 = {A_2}/{t}$.\\[0.1in]
{\bf Condition 1.} Unlike conditions 2 and 3, we are unable to satisfy condition 1~(in
either Proposition~\ref{prop3} or \ref{prop4})
without taking into account the particular representation of our group. In all cases we will need
to select an appropriate representation that will allow us to ensure that the conjugacy classes of
selected elements $\{\sigma_1, \ldots , \sigma_m\}$ of the Weyl group of a maximal torus $T$ give
rise to a strictly transitive set.  As noted before this representation must be minuscule. 
Using Proposition~\ref{local} and Corollary~\ref{cor1}, we can
guarantee that we can construct a differential equation $Y' = AY$ having singular points $p_1,
\ldots p_m$ so that the local formal Galois group at $p_i$ is isomorphic to the group generated by
$\sigma_i$ and $T$. Since the set $\{\overline{\sigma_1}, \ldots \overline{\sigma_m}\}$ will, by 
construction, be a strictly transitive set, condition 1 of Proposition~\ref{prop4} will be met.  In fact, for simply connected groups of type $\Al$ and $\Cl$
(\ie, $\SL_{\ell+1}$ and $\Sp_{2\ell}$) we will only need to use one element from the Weyl group
and the construction can be made so explicit that Proposition~\ref{prop3} can be applied.  We now
complete this argument for each of the above types of groups. 
\subsubsection{Proof of Proposition \ref{propsimply} for the Type $\boldsymbol{\Al}$} The simply connected group of type $\Al$ is $\SL_{\ell+1}$.  We
shall consider the usual representation of this group acting on a vector space of dimension
$\ell+1$ (in fact, all the fundamental representations are minuscule but we will only consider this one). In this representation the  diagonal 
elements form a maximal torus $T$ and  there are $\ell+1$ distinct weights.
The Weyl group is isomorphic to the group of unimodular permutation matrices, that is, to $\gS_{\ell+1}$,  and its action on the weights is the usual action of this permutation group on $\ell + 1$ elements. 
  In particular, there is an element $\sigma_1$ of the Weyl group that cyclically
permutes the roots. Clearly, $\{\overline{\sigma_1}\}$ forms a strictly transitive set. Using Proposition~\ref{local}, we can find an element $\overline{A}
\in \sl_{\ell+1}(C[t,t^{-1}])$ such that the Galois group of $Y' = \overline{A}Y$ is the group generated
by $\sigma_1$ and $T$.  Letting $B_3 = \overline{A}$, we see from Proposition~\ref{prop4} that together with
the matrices $B_1$ and $B_2$ already constructed, the set $\{B_1, B_2, B_3\}$ satisfies Proposition~\ref{propsimply}.  \\[0.1in]
The groups $\SL_{\ell+1}$ are particularly transparent but this is not the case of the other groups we will consider. For simply connected groups of the remaining types, the action of the Weyl group on weights for a given representation is best described in Lie-theoretic terms. We will in each case identify the Weyl group with the group generated by reflections associated to a set of simple roots in the Lie algebra. To prepare the reader for this
discussion we will show 
how one can find the element $\sigma_1$ above  using the Lie-theoretic description of the Weyl group.   
To do this we now
fix some notation.\\
[0.1in]
Let ${\g}$ be a semisimple Lie algebra of rank $\ell$ with Cartan subalgebra $\h$ and $\{\e_1, \ldots , \e_{\ell}\}$
a basis of the dual vector space $\h^*$. We will use an inner product on $\h^*$ for which the $\e_i$ form an
orthogonal basis and  write weights in terms 
of the $\e_i$ in order to do our calculations. \\[0.1in]
We will use the description of $\Al$ given on Planche I of Bourbaki (1968). The simple roots are given as
\[\alpha_i =   \e_i- \e_{i+1}, \ \ \ i = 1, \ldots , \ell \ \ . \]
We will only look at the minuscule representation $V(\obar_1)$ whose highest weight is $\obar_1 = \e_1 -
\frac{1}{\ell+1}\sum_{j=1}^{\ell+1}\e_j$.
This is the
standard representation and it is of dimension $\ell+1$.  We wish to now determine
which weights appear in this representation and how the Weyl group permutes them. We denote by $S_i$ the reflection
associated with the simple root $\alpha_i$, that is \[S_i(v) = v-\frac{2(v,\alpha_i)}{(\alpha_i,\alpha_i)}\alpha_i \ .\] 
Define $w_j = \e_j - \frac{1}{\ell+1}\sum_{j=1}^{\ell+1}\e_j$ and note that $w_1 = \obar_1$. We also
note that for all $i$, $1\le i \le \ell$, we have $(\alpha_i,\alpha_i) = 2$, so
\begin{eqnarray*}(1\leq i \leq \ell) \ \ \ S_i(w_j) & = & w_j - (w_j , \alpha_i) \alpha_i \\
& = & w_j-(\e_j - \frac{1}{\ell+1}\sum_{j=1}^{\ell+1}\e_j, \e_i -\e_{i+1})(\e_i-\e_{i+1})\\
    & = &
\begin{cases}
      w_j & \text{ if $j \neq i, i+1$}\\
      w_j -(\e_i - \e_{i+1}) & \text{ if $j = i$} \\
      w_j +(\e_i -\e_{i+1}) & \text{ if $j = i+1$}
      \end{cases}\\
& = &
\begin{cases}
      w_j & \text{ if $j \neq i, i+1$}\\
      w_{j+1} & \text{ if $j = i$} \\
      w_{j-1}& \text{ if $j = i+1$.}
      \end{cases}
  \end{eqnarray*}
Since $V(\omega_1)$ is minuscule, the set of weights appearing in this representation is precisely the orbit of $\obar_1$ under the Weyl group. 
The previous calculation shows that $\{w_i\}_{i=1}^{\ell+1}$ are among these and, comparing dimensions, we have that
this set is precisely  the set of weights. Again from this calculation, we see that $S_i =(i,i+1)$ where we
are thinking of $S_i$ as a permutation of the (subscripts of the) weights.  A calculation shows that
$S_1S_2\cdot\ldots\cdot S_\ell = (1,2,\ldots,\ell,\ell+1)$.  Let $\sigma_1$ be an element of the normalizer of a
maximal torus corresponding to this latter permutation (this can be found once one has elements corresponding to
the reflections $S_i$, see  Exercise 23.22 of Fulton and Harris, 1991). This
 is the element $\sigma_1$ introduced in the first paragraph.

\noindent We finish our discussion of $\SL_{\ell+1}$ by showing how 
we can make the  above construction  even more explicit.  
Prior to our work on Proposition~\ref{prop4}, the late A. Bolibrukh and 
R. Sch\"afke showed us how to explicitly construct differential equations with 
irreducible local Galois groups. Using their ideas in combination with the Weyl
 group approach, we may proceed as follows.
%
Let $A_{0.1} = (a_{i,j})$ be the matrix  defined by $a_{i+1, i} = 1$
for $i = 1,
\ldots ,\ell$ and $a_{i,j} = 0$ if $j+1 \neq i$ and let $A_{0,2}$ be the matrix with $1$ as the
$(1,\ell+1)$ entry and $0$ everywhere else. 
Note that  $A_{0,1}+A_{0,2}$ is a matrix whose eigenvalues are the $(\ell +1)^{st}$
roots of $1$. 
We now apply Corollary~\ref{cor1}. Select three points $p_0,p_1,
 p_3$  not in $\calS$,  with distinct projections on $\P1$. Let $A_1 = {\rm diag}(r_1,
\ldots, r_{\ell+1})$ where $r_1, \ldots , r_{\ell}$ are ${\bf Z}$-linearly independent {\it mod.} ${\bf
Z}$ and $r_{\ell+1} = -\sum_{j=1}^{\ell}$.  Let $A_2 = (a_{i,j})$ where $a_{i,j} = 0$ if $i \geq j$ and $a_{i,j} = 1$  if $i < j$.
Corollary~\ref{cor1} implies that one can produce an $A \in \calG(K)$ such that in terms of the
local coordinate $t$ at these points, the equation $Y' = AY$ has the following local expansions:
\begin{eqnarray}
\mbox{At }  p_0,\ \frac{dY}{dt} &=& (\frac{A_{0,1}}{t^2} + \frac{A_{0,2}}{t} +\mbox{ terms involving
}  t^j, \ j \geq 0)Y \label{e1}.\\
\mbox{At }  p_1,\ \frac{dY}{dt} &=& (\frac{A_1}{t} + \mbox{ terms involving }  t^j, \ j
\geq 0)Y \label{e2}. \\
\mbox{At }  p_2,\ \frac{dY}{dt} &=& (\frac{A_2}{t} + \mbox{ terms involving } t^j, \ j
\geq 0)Y\label{e3}. 
\end{eqnarray}
We now will check that the conditions of Proposition~\ref{prop3} hold. To see that there is a unique
slope at $p_0$,  let $g = \mbox{diag }(1,t^{1/(\ell + 1)}, t^{2/(\ell + 1)},
\ldots , t^{\ell/(\ell+1)})$. Note that for any matrix $(a_{i,j})$, we have that 
\[
g(a_{i,j})g^{-1} = (t^{\frac{i-j}{\ell+1}}a_{i,j}).
\]
Therefore 
\[
g[A] = gAg^{-1} + g'g^{-1} = \frac{A_{0.1}+A_{0,2}}{t^{2-\frac{1}{\ell+1}}} +
\mbox{ terms involving } t^j, \ j \geq
 2-\frac{1}{\ell+1}.\]
This is a  so-called  {\it shearing-transform} of the equation ${dY}/{dt}=AY$ at $p_0$. Since the matrix
$A_{0,1}+A_{0,2}$ is semisimple, there will be a unique
 slope, equal to  $2 - 1/(\ell+1)$ (Babbitt and  Varadarajan, 1983,  Proposition 4.2 
and the subsequent paragraphs).  Therefore,  as noted before, 
the equation will be irreducible over $C((t))$. \\[0.1in] 
Finally, at $p_1$ and $p_2$ the required conditions are obviously satisfied.  Therefore, the
equivariant equation $Y' = AY$ has Galois group $G^0$ over $K$ and so using the techniques of
Mitschi and Singer (2002) or Hartmann (2002) one can construct an equation having Galois group $H \ \semi G^0$
over $C(x)$. In particular, the example given in Section 2 was constructed in the above manner and
so has Galois group $\SL_2$ over $K = C(x, \sqrt{x})$ (another example with this group is given {\it via}
an {\em ad hoc} construction 
by Hartmann (2002), p.~42, and in Section~\ref{alternate})\\[0.2in]
\subsubsection{Proof for  $\boldsymbol{\Cl}$} The simply connected groups of this type are the groups ${\Sp_{2\ell}}$ (see the tables for ($C_l$) given by Tits, 1967, p. 32).
A set of
simple roots of $\Cl$ are
\begin{eqnarray*}
\alpha_i & = &\e_i-\e_{i+1},  \ i = 1, \ldots , \ell-1,\\
\alpha_{\ell} & = & 2\e_\ell
\end{eqnarray*}
(see Planche III of Bourbaki, 1968).  
The only minuscule weight is $\obar_1 = \e_1$, corresponding to the standard representation $V(\obar_1)$
 of $\Sp_{2\ell}$
 which has dimension $2\ell$. We again denote by $S_i$ the reflection across the simple root $\alpha_i$
and will calculate the action of these reflections on the $\e_i$. Noting that for $1\leq i \leq\ell-1, (\alpha_i,
\alpha_i) = 2$ and $(\alpha_\ell, \alpha_\ell) = 4$ we have
\begin{eqnarray*}(1\leq i \leq \ell - 1) \ \ \ S_i(\e_j) & = &\e_j - (\e_j,  \e_i - \e_{i+1}) (\e_i-\e_{i+1}) \\
    & = &\begin{cases}
      \e_j & \text{ if $j \neq i, i+1$}\\
      \e_{i+1} & \text{ if $j = 1$} \\
      \e_i & \text{ if $j = i+1$}
      \end{cases}\\
      & & \\
   S_\ell(\e_j) & = & \e_j - \frac{2(\e_j ,2\e_\ell)}{4} 2
\e_\ell \\
    & = &\begin{cases}
      \e_j & \text{ if $j \neq \ell$}\\
     -\e_{\ell} & \text{ if $j = \ell$.} 
      \end{cases}
  \end{eqnarray*}
Since $V(\obar_1)$ is minuscule, the set of weights appearing is precisely the orbit of $\obar_1$ under the Weyl group. 
The above calculation shows that this orbit contains $\{\pm \e_i\}_{i = 1}^{\ell}$ and so by comparing dimensions we
see that it is precisely the set of weights. Let us give the labels  $\e_1 = 1, \ldots , 
\e_{\ell} = \ell, -\e_1 = \ell+1, \ldots, -\e_{\ell} = 2\ell$.   Rewriting the $S_i$ as permutations of these weights, the
above calculation shows that $S_1\cdot\ldots\cdot S_{\ell-1} = (1,2,\ldots , \ell)(\ell+1, \ldots ,2\ell)$ and
$S_1\cdot\ldots\cdot S_{\ell} = (1,2,\ldots , 2\ell)$.  Let $\sigma_1$ be an element of the normalizer of a maximal
torus of $\Sp_{2\ell}$ that yields this latter permutation. 
We see that
$\sigma_1$ acts transitively on the set of weights and so $\{\overline{\sigma_1}\}$ forms a strictly transitive set.
One now proceeds as in the first paragraph of the discussion  of $\Al$.\\[0.1in]
We note that we can make this as explicit as the example in the discussion of $\Al$.  To do this we let $U =
(a_{i,j})$ be the marix defined by $a_{i+1, i} = 1$ for $i = 1, \ldots , \ell-1$ an $a_{i,j} = 0$ if $j+1 \neq i$.
and let $V$ be the matrix with $1$ as the $(1,\ell)$ entry and $0$ everywhere else. Let
\[A_{0,1} = \left(\begin{array}{cc} U&0\\V&-U\end{array}\right) \ \ \ \ 
A_{0,2} = \left(\begin{array}{cc} 0&(-1)^\ell V\\0&0\end{array}\right).\]
Note that  each of these matrices is in $\Sp_{2\ell}$ and $A_{0,1}+A_{0,2}$ is a matrix whose eigenvalues are the
$2\ell^{th}$ roots of unity.  We also define $A_1 = {\rm diag}({\rm diag}(r_1, \ldots , r_\ell),-{\rm diag}(r_1, \ldots ,
r_\ell))$ where $r_1, \ldots , r_{\ell-1} \in C$ are $\bf Z$-linearly independent {\it mod.} $\bf Z$ and $r_n = -\sum r_i$.
 Finally, we let $A_2$ be any nilpotent matrix in $\Sp_{2\ell}$ such that $\Ad(u)$ has an $r$-dimensional eigenspace
corresponding to $1$.  One then proceeds as in the case of $\Al$.
\subsubsection{Proof for  $\boldsymbol{\Dl}$} The simply connected groups of this type are the groups ${\Spin_{2\ell}}$.
These are double covers
of the groups $\SO_{2\ell}$. We  claim that it suffices to prove the analog of Proposition \ref{propsimply}  for all (non simply connected) groups $\SO_{2\ell}$.   To see this, assume 
that we have verified the result of Proposition \ref{propsimply} for one of these groups, say $G$.  Let $\pi:G' \rightarrow G$ 
 be a simply connected
covering  with $d\pi:\calG' \rightarrow \calG$ the associated map of Lie algebras and $\rho': G' \rightarrow \GL_{N'}$ be a faithful representation.  The map $d\pi$ will be an isomorphism from $\calG'$ onto $\calG$. Let $B'_i = d\pi^{-1}(B_i)$.  We claim that 
the $B'_i$ satisfy the conclusions of Proposition \ref{propsimply}.  Let $B' \in \gl_{N'}(K)$ satisfy  hypotheses 1 and 2 of the proposition and let $K(g)$ be a Picard-Vessiot
extension for $Y' = B'Y$ with $g \in G'$. Then $\pi(g)$ satisfies $Y' = d\pi(B)Y$ and  this latter equation has Galois group
$G$ by the above proposition.  Since the only subgroup of $G'$ mapping onto $G$ {\it via} $\pi$ is $G'$, Proposition 5.3 of Mitschi and Singer (2002)
implies that the Galois group of $K(g)$ over $K$ is $G'$.\\[0.1in]
\noindent A set of simple roots in the case of $\SO_{2\ell}$ is  
\[\alpha_1  = \e_1-\e_2, \ldots , 
\alpha_{\ell-1} = \e_{\ell-1}-\e_\ell, \alpha_{\ell} = \e_{\ell-1} + \e_{\ell}\]
(see Planche IV of Bourbaki, 1968) 
 We shall
only consider the representation $V(\obar_1)$ corresponding to $\obar_1 = \e_1$.  This is the standard
representation of $\SO_{2\ell}$ of dimension $2\ell$ (and is not a faithful representation of $\Spin_{2\ell}$). 
Noting that $(\alpha_i, \alpha_i) = 2$ for $i = 1, \ldots \ell$, we 
again make the calculation
\begin{eqnarray*}(1\leq i \leq \ell - 1) \ \ \ S_i(\e_j) & = & \e_j - (\e_j , \e_i-\e_{i+1}) (\e_i - \e_{i+1}) \\
    & = &\begin{cases}
      \e_j & \text{ if $j \neq i, i+1$}\\
      \e_{i+1} & \text{ if $j = 1$} \\
      \e_i & \text{ if $j = i+1$}
      \end{cases}\\
      & & \\
   S_\ell(\e_j) & = & \e_j - (\e_j ,\e_{\ell-1} + \e_{\ell}) (\e_{\ell-1} + \e_{\ell}) \\
    & = &\begin{cases}
      \e_j & \text{ if $j \neq \ell, \ell-1$}\\
     -\e_{\ell} & \text{ if $j = \ell-1$} \\
      \e_{\ell-1} & \text{ if $j = \ell$.}
      \end{cases}
  \end{eqnarray*}
Once again, we can conclude that the orbit of $\obar_1$ is $\{\pm \e_i\}_{i=1}^\ell$.
  If we give the labels $\e_1 = 1, \ldots , \e_{\ell} = \ell, -\e_1 = \ell+1, \ldots, -\e_{\ell} = 2\ell$, then the $S_i$
  correspond to the following permutations:
  \begin{eqnarray*}
  S_i & = & (i, i+1)(\ell+i,\ell+i+1), \ \ \ 1\leq i \leq \ell-1\\
  S_\ell & = & (\ell-1,2\ell)(\ell, 2\ell - 1).
  \end{eqnarray*}
  From this one can show
   \begin{eqnarray*}
  S_1S_2\cdot\ldots\cdot S_{\ell-1} & = & (1,\ldots ,\ell)(\ell+1,\dots ,2\ell) \\ 
  S_1S_2\cdot\ldots\cdot S_{\ell} \ \ \ & = & (1,\ldots ,\ell-1,\ell+1,\dots ,2\ell-1)(\ell, 2\ell).
  \end{eqnarray*}
Let $\sigma_1$ and $\sigma_2$ be  elements of the normalizer of a maximal
torus $T$ of $\SO_{2\ell}$ that yield these permutations.
Then $\{\overline{\sigma_1}, \overline{\sigma_2}\}$ is a strictly transitive set.\\[0.1in]
We  construct $B_3$,  $B_4$ for the groups generated by $\sigma_1$ and $T$, by  $\sigma_2$ and $T$ respectively.  
We then have that  $B_1, B_2,B_3,B_4$ yield the conclusion of Proposition~\ref{propsimply}.  

%

%
\subsubsection{Proof for  $\boldsymbol{\esi}$} The positive roots and a choice of simple roots are given on Planche V of Bourbaki (1968)
but it is slightly easier to make the computation using the form given by Fulton and Harris (1991,  on
p. 332-333),  (where $L_i$ is used for $\e_i$ and $\omega$ is used for $\obar$).
 There are two possible minuscule representations; those having highest weight $\obar_1 = \frac{2\sqrt{3}}{3}L_6$ and 
$\obar_6 = L_5 + \frac{\sqrt{3}}{3}$ ({\it loc. cit. }, formulas for ($E_6$)  p. 333 and p. 528).  We will consider the $27$-dimensional representation 
corresponding to $\omega_1$. 
According to the tables for ($E_7$) given by  Tits (1967, p. 47), the representation associated to $\omega_1$ is a faithful representation 
of the simply connected group of this type.\\[0.1in] 
Using a Maple package developed by the first author (see the link {\tt www.math.ncsu.edu/
$\sim$singer/papers/weyl\_permutation.html} to download the software and
reproduce the calculation), we calculated the $27$ weights associated with this representation and calculated the permutations induced by the reflections
given by the simple roots. We were able to determine the cycle structure of all elements in the permutation group generated by  these elements and 
found that there were elements with cycle structure $[12,12,3]$ (that is, a product of two $12$-cycles and one $3$-cycle)
and $[9,9,9]$. A simple calculation shows that the associated conjugacy classes
act strictly transitively. We refer to the above web page for details of the computation. Letting $\sigma_1$ and $\sigma_2$ be elements of the normalizer
of a maximal torus having these cycle structures, we proceed as above to find the $B_i$.
\subsubsection{Proof for $\boldsymbol{\ese}$} We again use the positive roots and 
a choice of simple roots as given by Fulton and Harris (1991, on p.~333). 
There is only one possible
minuscule representation: the representation having highest weight $\omega_7 = L_6 + \frac{\sqrt{2}}{2}L_7$. This representation has dimension
$56$ and is a faithful representation of the associated simply connected group, as appears in the tables p. 47 of Tits (1967). 
A calculation similar to the one above shows that there are two elements of the Weyl group whose permutation structure is given
by $[18,18,18,2]$ and $[14,14,14]$ respectively.  The associated conjugacy classes are strictly transitive. We again refer to the web page for 
the details.
\subsubsection{Proof for Other Types} Groups of type ${\rm G}_2$,  ${\rm F}_4$ and ${\rm E}_8$ do not have minuscule representations and so the above methods do not apply.
The spin  representation (corresponding to highest weight $\omega_\ell$) is a minuscule representation for groups of 
type $\Bl$. We have calculated the action of the Weyl group on the weights of this representation and are able to produce strictly transitive
sets of permutation conjugacy classes for ${\rm B}_2$, ${\rm B}_3$, ${\rm B}_5$ and ${\rm B}_7$ and can show that such sets do not exist for ${\rm B}_4$. We do not have 
definitive results for general groups of type $\Bl$. 
 
\subsection{Equivariant Equations for $\boldsymbol{G^0 =  G_1 \cdot\ldots\cdot G_r}$   where each $\boldsymbol{G_i}$ is  of Type 
 $\boldsymbol{\Al}$, $\boldsymbol{\Cl}$, $\boldsymbol{\Dl}$,
 $\boldsymbol {\esi}$ or
$\boldsymbol{\ese}$}
At the beginning of this section, we showed how one can reduce the problem under consideration to finding equivariant equations for {\em simply connected} groups of the same type,
that is groups $G^0 = \prod G_i$ where each $G_i$ is simply connected and of type $\Al, \Cl, \Dl, {\rm E}_6$ or ${\rm E}_7$. We shall restrict ourselves
to groups of this form. Let $\calG_i$ denote the Lie algebra of $G_i$ and 
$\calG =
\oplus_{i = 1}^r
\calG_i$. For each $i$, let $B_{i,1} , \ldots ,B_{i,m_i}$ be the elements guaranteed to exist by Proposition~\ref{propsimply}.
Let $\{p_{i,1}, \ldots , p_{i,m_i}\}_{i=1}^r$ be points on $\curve\backslash \calS$ having
distinct projections. Corollary~\ref{cor1} implies
that  one can  find an equivariant $B = {\rm diag}(B_1, \ldots B_r) \in \calG(K) = \oplus \calG_i(K)$ with  $B_i \in \calG_i(K)$ such that in terms of the local coordinate $t$ at the
point $p_{i,j}$, the equation $Y' =B_iY$ has the form ${dY}/{dt} = (B_{i,j} + \mbox{ higher order terms})Y$
and such that the equation is non-singular at the points $p_{j,k}$ for $j \neq i$. \\[0.1in]
Let $E$ be the Picard-Vessiot extension of $K$ corresponding to $Y' = BY$.  Since $B \in \calG(K)$,
the proof of Proposition 1.31 of van der Put and Singer (2003) shows that we can assume that $K$ is generated by
the entries of an element $g \in G^0(E)$ such that $g' = Bg$. Writing $g = (g_1, \ldots, g_m)$ where
each $g_i $ is in  $G_i$, we have that  $g_i' = B_ig_i$ and so $E$ contains the Picard-Vessiot extension
$E_i=K(g_i)$ of $K$ corresponding to each of the equations. From Proposition~\ref{propsimply}, 
we know that the Galois group  of $Y' = B_iY$ over $K$ is $G_i$. We shall now show
that the Galois group $G'$ of $Y' = BY$ over $K$ is $G^0$.\\
[0.1in]
Since $A \in \calG(K)$, we have that $G' \subset G$ (see Proposition 1.31 of van der Put and Singer, 2003).  Assume
that $G'
\neq G$.  We will show that this implies that there exist indices $i \neq j$ such that $E_i$ lies in
an algebraic extension of $E_j$.  We will see that  comparing the local behavior of solutions of the
corresponding differential equations at some $p_{i,k}$ will yield a contradiction.\\
[0.1in]
 A result of Kolchin (1968) (see also Exercise 8, Chapter V.23 of
Kolchin, 1973) implies that there are indices  $i\neq j$ and a homomorphism (defined over $C$)  $f:G_i
\rightarrow G_j/Z(G_j)$, where $Z(G_j)$ is the center of $G_j$, such that for every $h=(h_1, \ldots
,h_m) \in G', \ f(h_i) =
\pi(h_j)$, where $\pi$ is the canonical homomorphism $G_j \rightarrow G_j/Z(G_j)$.  Note that since
$G_i$ and $G_j$ are simple, the kernels of $f$ and $\pi$ are finite. \\
[0.2in]
We now apply the maps $f$ and $\pi$ to the element $g=(g_1, \ldots,g_m) \in G^0(E)$ defined above.
Since $f(g_i)=\pi(g_j)$, we have that $E_i$ and $E_j$ share the common subfield
$K(f(g_i))=K(\pi(g_j))$. Furthermore, $E_i$ and $E_j$ are algebraic extensions of this field since
the kernels of $f$ and $\pi$ are finite. Therefore $E_i$ is contained in an algebraic extension of
$E_j$. \\
[0.2in]
By construction, $Y' = B_jY$ is non-singular at each of the $p_{i,k}$ and so  the solutions of $Y' = B_jY$ at
$p_{i,k}$ have components in $C((t))$ where $t$ is the local parameter at $p_{i,k}$. Therefore we can embedd $E_j$
into $C((t))$.  This implies that $Y' = B_iY$ has a fundamental set of solutions in an algebraic
extension of $C((t))$ and so must be regular singular at this point (see Exercise 3.29 of van der Put and Singer, 2003). By construction, one of the points $p_{i,k}$ is not a regular singular point and so this is a contradiction. Therefore
the Galois group of $Y' = BY$ is $G$.

\section{An Alternate Construction for Finite Extensions of $\boldsymbol{\SL_2}$}\label{alternate} In this section
we present
an alternate method for constructing linear differential equations whose Galois groups are finite
extensions of $\SL_2$. In the previous sections, we considered groups of the form $H \ \semi G^0$,
$H$  a finite group and $G^0$ of the type considered above, and showed that for any realization of
$H$ as a Galois group of an extension $K$ of $C(x)$, we could find an equivariant $A$ such that $Y'
= AY$ had Galois group $G^0$ over $K$.  The construction described here begins by constructing a
suitable $K$ and so does not work over any such $K$.  On the other hand, it introduces fewer
singularities and uses group theoretic facts that may be of independent interest.  This construction
was motivated by the Example given by Julia Hartmann (2002, p. 42).\\
[0.2in]
We begin with a modification of Lemma 5.11 of Borel and Serre (1964-65) (see also Lemma 10.10 of Wehrfritz, 1973).  For any algebraic group $G$ we
define $\Int:G\rightarrow
\Aut(G^0)$ to be the map that sends an element to the automorphism resulting from conjugation by
that element.
\begin{lem} \label{lem5.1} Let $G$ be a linear algebraic group, $B$ a Borel subgroup of $G$ and $T$
a maximal torus of $B$. \\
[0.1in] 1.~There exists a finite subgroup $W$ of $G$  such that 
$W$ normalizes $B$ and $T$ and  the natural projection $W \rightarrow G/G^0$ is surjective. \\
[0.1in]
2. If, in addition, $G^0$ is semisimple and all automorphisms of $G^0$ are inner, then
 $\Int(W) \subset \Int(T)$ and so $\Int(W)$ is a finite abelian group. 
If $G^0 =
\SL_2$ or $\PSL_2$, then $\Int(W)$ is cyclic.

\end{lem}
\begin{proof} 1.~Let $N_G(B)$ be the normalizer of $B$ in $G$ and $N_G(B,T)$ be the subgroup of
elements of $N_G(B)$ that normalize $T$ as well. Since all Borel subgroups of $G$ lie in $G^0$ and
are conjugate in $G^0$ (see Theorem 21.3 of Humphreys, 1975), we have that for any $g \in G$ there exists
an $h \in G^0$ such that $gBg^{-1} = hBh^{-1}$. Therefore $h^{-1}g \in N_G(B)$ and we can conclude
that $G = N_G(B) \cdot G^0$. 
Using the fact that the maximal tori of $B$ are all conjugate in $B$ ({\it loc. cit.} Theorem 19.3), we also have  that $N_G(B) = N_G(B,T) \cdot B$.  Lemma 10.10 of 
Wehrfritz (1973) implies that there exists a finite subgroup $W$ of $N_G(B,T)$ such that the natural
projection $W \rightarrow N_G(B,T)/N_G(B,T)^0$ is surjective.  We then have that the projection
$W\rightarrow G/G^0$ is surjective as well.\\
[0.2in]
2. We refer to the proof  of  Theorem 27.4 of Humphreys (1975).  Since all automorphisms of $G^0$ are
inner, for any element $w \in W$ there is an element $h\in G^0$ such that for all $g \in G^0,
wgw^{-1} = hgh^{-1}$.  Since $w$ normalizes $B
\mbox{ and } T$, we have that  $\Int(W)\subset\Int(N_{G^0}(B,T))$. Since $B$ is a Borel subgroup, we
have that $N_{G^0}(B,T) \subset N_{G^0}(B) = B$. An element of $B$ that normalizes $T$ must lie in
$T$  ({\it loc. cit.} Proposition 19.4, Corollary 26.2A) so $\Int(W) \subset
\Int(T)$.  The final statement follows from the fact that a maximal torus of these groups has
dimension $1$.\end{proof}

Let $G$ be a linear algebraic group with $G^0$ semisimple and let $\calG$ be its Lie algebra. If $T$
is a maximal torus of $G^0$, then its Lie algebra $\calT$ is a Cartan subalgebra of $\calG$ and we
can decompose
\[
\calG = \calT \oplus \prod_{\alpha \in \Phi} \calG_{\alpha}
\]
where $\Phi$ are the roots of $\calG$ which we consider as multiplicative characters on $T$. If $W$
is the finite group described in Corollary~\ref{lem5.1}, then for any $\alpha \in \Phi$ and $w \in
W$ we define $\alpha(w) =
\alpha(t)$ for any $t \in T$ such that $\Int(w) = \Int(t)$. Since the elements of $\Phi$ factor
through $\Int:G
\rightarrow \Aut(G)$, each root in this way defines a multiplicative character on $W$.

\begin{ex} {\em Let $G^0 = \SL_2$ and assume that $T$ is the subgroup of diagonal matrices. As usual
we let
\[
h = \left(
\begin{array}{cc} 1 & 0 \\
0& -1 
\end{array} \right), \ \ e =
\left(
\begin{array}{cc} 0 & 1 \\
0& 0
\end{array} \right), \ \ f = \left(
\begin{array}{cc} 0 & 0 \\
1& 0
\end{array} \right) \ .   
\]
$\calT$ is spanned by $h$ and there are two roots $\alpha$ and $-\alpha$ with $\calG_{\alpha}$ being
spanned by $e$ and $\calG_{-\alpha}$  by $f$. Furthermore, considering the roots as characters on
$T$, we have that
\[
\alpha \left(
\begin{array}{cc} a & 0 \\
0& a^{-1} 
\end{array} \right) = a^2, \ \ \ \ -\alpha
\left(
\begin{array}{cc} a & 0 \\
0& a^{-1} 
\end{array} \right) = a^{-2} 
\]
Let  $G = \SL_2 \ \semi \{ 1, -1\} $ where the action of $-1$ on $\SL_2$ is given by conjugation by
the matrix
\[
\left(
\begin{array}{cc} \sqrt{-1} & 0 \\
0& -\sqrt{-1} 
\end{array} \right) \ .
\]
We then have that $\Int(W)= \Int(H)$ where $W=\{1,-1\}$ and $H$ is the order four cyclic subgroup of
$\SL_2$ generated by the above matrix. Note that 
$\alpha$ can be considered as the character on $W$ given by $\alpha(-1) = -1$. \QED}
\end{ex}
Let $G$ be an algebraic group with $G^0 = \SL_2$ and let $W$ be as in Lemma~\ref{lem5.1}. We shall
construct a differential equation having $W\ \semi G^0$ as its Galois group over $C(x)$. The group
$G$ will then be the Galois group of a subfield $\tilde{E}$ of the Picard-Vessiot extension $E$ of
$C(x)$ corresponding to this former equation.  \\
[0.2in]
We now use the notation $G$ to denote the group $\SL_2 \ \semi W$ and $\calG$ to denote the Lie
algebra of $G$.  Conjugation by an element of $W$ induces an automorphism of $G$ and also an
automorphism of $\calG$, which we again denote by conjugation.  If $K$ is any field containing $C$
and $X = ah+be+cf \in
\calG (K), \ a,b,c \in K$, then, for $w \in W$
\[
w^{-1}Xw = ah + b\alpha(w^{-1})e + c(-\alpha(w^{-1})f)
\]
Note that if the image of $W$ in $\Aut(G)$ has order $n$, then $\alpha$ maps $W$ onto the group of
${ n}^{ th}$ roots of unity.  We identfy this with the Galois group of $C(x, x^{1/n})$.
 Let $K$ be a Galois extension of $C(x)$ with Galois group $W$ such that the fixed field of the
kernel of $\alpha$ is $C(x, x^{1/n}), x' = 1$ and the action of $W$ on this latter field is given by
$\alpha$.  Theorem 7.13 of Volklein (1996) implies that such a field exists.   Let
\[
\tilde{A} = x^{-1/n}e \ + \ x^{1/n} f \ + \ x^2 h \ .
\]
To construct a differential equation whose Galois group is $G$, Proposition 5.2 of Mitschi and Singer (2002)
implies that it is enough to prove the following proposition.
\begin{prop}
 $\tilde{A}$ is equivariant and the differential Galois group of 
 $Y' = \tilde{A}Y$ over $C(x, x^{1/n})$ is $\SL_2$.
\end{prop}
\begin{proof}  To prove the claim about the Galois group, we make a change of variables $x=z^n$.  We
then get a new equation $\frac{dY}{dz} = AY$ where
\[
A = n(z^{n-2}e + z^n f + z^{3n-1} h)
\]
We will use the techniques of Mitschi and Singer (1996) to show 
that $\frac{dY}{dz} = AY$ has differential Galois group $\SL_2$ over $C(z)$.\\
[0.2in]
Assuming that this latter fact is true, we claim that the differential 
Galois group of  $Y' = \tilde{A}Y$ over $C(x, x^{1/n})$ is $\SL_2$. To see this, let $K$ be a
Picard-Vessiot 
extension of $C(z) = C(x^{1/n})$  for $\frac{dY}{dz} = AY$. 
Since $K$ has no new $\frac{d}{dz}$-constants, it has no new 
$\frac{d}{dx}$-constants. Furthermore, since the elements of $\SL_2$ commute with $\frac{d}{dz}$ and
leave $C(z)$ fixed,
 they will commute with $\frac{d}{dx} = \frac{1}{nz^{n-1}}\frac{d}{dz}$. Therefore, $\SL_2$ is a
subgroup of the differential Galois group $G$ of $K/C(z)$ with respect to $\frac{d}{dx}$. From the
form of $\tilde{A}$, we see that $G \subset \SL_2$, so the claim is proved.\\
[0.2in]
We now proceed to show that $\frac{dY}{dz} = AY$ has differential Galois group $\SL_2$ over $C(z)$.
Since $C^2$ is a Chevalley module for $\SL_2$, Lemma 3.3 of Mitschi and Singer (1996) implies that it is enough
to show that if
\[
c = c_{3n-1}z^{3n-1} + \ldots + c_0 \in C[z] \ {\rm and } \ w = w_mz^m + \ldots + w_0 \in C^2
\otimes C[z]
\]
and 
\begin{eqnarray}
\label{eqn1} w' - n[z^{n-2}e + z^n f + z^{3n-1} h - cI]w& = &0
\end{eqnarray}
then $w = 0$.  \\
[0.2in]
To simplify notation, we let  $ u =
\left(
\begin{array}{c} 1\\
0 
\end{array}\right)$ and
 $ v = \left(
\begin{array}{c} 0\\
1 
\end{array}\right)$. These generate the root spaces of $\SL_2$ and $eu =0, ev=u, fu=v, 
fv=0$,  and t $u,v$ are eigenvectors of $h$.  The proof of the above fact proceeds by considering
the coefficients of powers of $z$ in equation~(\ref{eqn1}).  The highest power of $z$ that can
appear is $z^{3n-1+m}$ and its coefficient is
\[
n(c_{3n-1}I-h)w_m =0 \ .
\]
We therefore have that $w_m$ is an eigenvector of $h$ and so we can assume that $w_m = u, c_{3n-1} =
1$ or $w_m = v, c_{3n-1} = -1$.  Let us assume that $w_m = u$ and $c_{3n-1} = 1$.  We shall write $w
=pu+qv$ where $p,q \in C[z], \ p = z^m +\mbox{ lower degree terms, and } q =
\mbox{ a polynomial of degree at most } m-1$. Substituting $w = pu+qv$ into equation (\ref{eqn1}), 
we have:
\[
(p' - nz^{n-2}q - nz^{3n-1}p + ncp)u + (q' - nz^np +nz^{3n-1}q +nqc)v = 0
\]
and therefore
\begin{eqnarray}
\label{eqn2} p' - n[z^{3n-1} -c]p & = &nz^{n-2}q\\
\label{eqn3}q' + n[z^{3n-1} +c]q& = & nz^np
\end{eqnarray}
The right hand side of equation~(\ref{eqn3}) has degree $n+m$. Since $z^{3n-1} + c$ has degree 
$3n-1$, $q$ must have degree  $m-2n+1$. Therefore the right hand side of equation~(\ref{eqn2}) has
degree $m-n-1$, while the degree of $[z^{2n-1} -c]p$ is at least $m$ if $z^{3n-1} - c \neq 0$. 
Therefore we have $c = z^{3n-1}$ and  that $p' = nz^{n-2}q$. Comparing degrees in this last
equation, we have $m-1 = n-2 + m-2n+1 = m-n-1$ so $n = 0$, a contradiction, unless $w=0$. If $w_m =
v$ and $c_{3m-1} = -1$ one argues in a similar way to also show that $w=0$. 
\end{proof}

\bigskip
\noindent
We wish to thank the late A.  Bolibrukh and R. Sch\"afke for showing us how to select the local form
of a differential equation so that Katz's Criterion ensures it is irreducible.

\medskip
\noindent
The third author would like to thank the Institut de Recherche Math\'ematique
Avanc\'ee, Universit\'e Louis Pasteur et C.N.R.S., for its hospitality and support during the
preparation of this paper. The preparation was also supported by NSF Grant CCR- 0096842.

\vspace{.4in}
\centerline{{\bf REFERENCES}}

\vspace{.3in}
\noindent
 Babbitt, D.~G.,  Varadarajan, V.~S. (1983). Formal reduction of meromorphic differential equations:a group theoretic view. {\em Pacific J. Math.}, 109(1):1--80.

\bigskip
\noindent
Bore, A.,  Serre, J.-P. (1964).
Th\'eor\`emes de finitude en cohomologie galoisienne.
{\em Comment. Math. Helv.}, 39:111--164.

\bigskip
\noindent
Bourbaki, N. (1968).
{\em {Groupes et Alg\`ebres de Lie, Chaps. 4,5 and 6}}.
Masson, Paris.

\bigskip
\noindent
Bourbaki, N. (1990).
{\em {Groupes et Alg\`ebres de Lie, Chaps. 7 and 8}}.
Masson, Paris.

\bigskip
\noindent
Fulton, W., Harris, J. (1991).
{\em {Representation theory. A first course}}, volume 129 of {\em
  Graduate Texts in Mathematics}.
Springer, New York1.

\bigskip
\noindent
Hartmann, J. (2002).
{\em On the inverse problem in differential {Galois} theory}.
Thesis, University of Heidelberg.
(available at {\tt www.ub.uni-heidelberg.de/archiv/3085}).

\bigskip
\noindent
Hochschild, G. (1976).
{\em {Basic Theory of Algebraic Groups and Lie ALgebras}}, volume~75
  of {\em Graduate Texts in Mathematics}.
Springer Verlag, Berlin.

\bigskip
\noindent
Hrushovski, E. (2002).
Computing the {Galois} group of a linear differential equation.
In {\em {Differential Galois Theory}}, volume~58 of {\em Banach
  Center Publications}, pages 97--138. Institute of Mathematics, Polish Academy
  of Sciences, Warszawa.

\bigskip
\noindent
Humphreys, J. (1975).
{\em {Linear Algebraic Groups}}.
Graduate Texts in Mathematics. Springer-Verlag, New York.

\bigskip
\noindent
Katz, N. (1987).
On the calculation of some differential {Galois} groups.
{\em Inventiones Mathematicae}, 87:13--61.

\bigskip
\noindent
Kolchin, E. R. (1968).
Algebraic groups and algebraic dependence.
{\em American Journal of Mathematics}, 90:1151--1164.

\bigskip
\noindent
Kolchin, E. R. (1973).
{\em {Differential Algebra and Algebraic Groups}}.
Academic Press, New York.

\bigskip
\noindent
Kovacic, J. (1969).
The inverse problem in the {Galois} theory of differential fields.
{\em Annals of Mathematics}, 89:583--608.

\bigskip
\noindent
Kovacic, J. (1971).
On the inverse problem in the {Galois} theory of differential fields.
{\em Annals of Mathematics}, 93:269--284.

\bigskip
\noindent
Lang, S. (1993).
{\em {Algebra}}.
Addison Wesley, New York, 3rd edition.

\bigskip
\noindent
Mitschi, C., Singer, M. F. (1996).
\newblock Connected linear groups as differential {Galois} groups.
\newblock {\em Journal of Algebra}, 184:333--361.

\bigskip
\noindent
Mitschi, C., Singer, M. F. (2002).
Solvable-by-finite groups as differential {Galois} groups.
{\em Ann. Fac. Sci. Toulouse}, XI(3):403 -- 423.

\bigskip
\noindent
van~der Put, M., Singer, M. F. (1997).
{\em {Galois Theory of Difference Equations}}, volume 1666 of {\em
  Lecture Notes in Mathematics}.
Springer-Verlag, Heidelberg.

\bigskip
\noindent
van~der Put, M., Singer, M. F. (2003).
{\em {Galois Theory of Linear Differential Equations}}, volume 328 of
  {\em Grundlehren der mathematischen Wissenshaften}.
Springer, Heidelberg.

\bigskip
\noindent
Singer, M. F. (1993).
\newblock Moduli of linear differential equations on the {Riemann} sphere with
  fixed {Galois} group.
\newblock {\em Pacific Journal of Mathematics}, 106(2):343--395.

\bigskip
\noindent
Tits, J. (1967).
{\em {Tabellen zu den einfachen Lie Gruppen und ihren
  Darstellungen}}, volume~40 of {\em Lecture Notes in Mathematics}.
Springer-Verlag, Heidelberg.

\bigskip
\noindent
Volklein, H. (1996).
\newblock {\em {Groups as Galois Groups}}, volume~53 of {\em Cambridge Studies
  in Advanced Mathematics}.
\newblock Cambridge University Press, Cambridge.

\bigskip
\noindent
Wehrfritz, B. A. F. (1973).
{\em Infinite Linear Groups}.
Ergebnisse der Mathematik. Springer-Verlag, Berlin.

\end{document}